\documentclass[11pt,a4paper] {article}
\usepackage{amssymb}
\usepackage{amsmath}
\usepackage{amsthm}
\usepackage{dsfont}
\usepackage{graphicx,epsfig}
\usepackage{float}
\usepackage{color}
\usepackage{capt-of}
\usepackage{tikz}
\usepackage[normalem]{ulem}

\newcommand{\R}{\mathbb R}

\newcommand{\EE}{\mathbb E}
\newcommand{\PP}{\mathbb P}

\newcommand{\cP}{\mathcal P}

\newcommand{\ds}{\displaystyle}

\newcommand{\cU}{\mathcal U}
\newcommand{\cF}{\mathcal F}
\newcommand{\cH}{\mathcal H}
\newcommand{\dt}{\Delta t}
\newcommand{\dk}{\Delta k}
\newcommand{\dz}{\Delta z}

\DeclareMathOperator*{\argmax}{argmax}

 \newtheorem{theorem}{\textbf{Theorem}}[section]
 \newtheorem{remark}[theorem]{\textbf{Remark}}

 \newtheorem{proposition}[theorem]{\textbf{Proposition}}
 \newtheorem{acknowledgement}{\textbf{Acknowledgement}}
 \newtheorem{assumption}[theorem]{\textbf{Assumption}}

 \definecolor{mypurple}{RGB}{140,0,255}
\definecolor{myred}{rgb}{255,0,0}

\definecolor{mydarkturquoise}{RGB}{0,206,209}
\definecolor{mydeeppink}{RGB}{255,20,147}
\definecolor{darkblue}{RGB}{0,0,140}
\definecolor{blue2}{RGB}{0,0,0}
\definecolor{middleblue}{RGB}{0,0,71}
\definecolor{light-gray}{gray}{0.9}

\definecolor{ProcessBlue}{cmyk}{1,0,0,0.40}
\definecolor{Black}{cmyk}{0,0,0,1}
\definecolor{Red}{cmyk}{0,1,1,0.2}
\definecolor{Green}{cmyk}{0.9,0,1,0}
\definecolor{Orange}{cmyk}{0,0.61,0.87,0.5}
\definecolor{Fuchsia}{cmyk}{0.47,0.91,0,0.06}
\definecolor{PineGreen}{cmyk}{0.92,0,0.59,0.30}

\title{A class of short-term models for the oil industry addressing speculative storage}
\author{Yves Achdou  \thanks { Universit{\'e} de Paris and  Sorbonne Universit{\'e}, CNRS, Laboratoire Jacques-Louis Lions, (LJLL), F-75006 Paris, France, achdou@ljll-univ-paris-diderot.fr}
,  Charles Bertucci \thanks{CMAP, {\'E}cole polytechnique. charles.bertucci@polytechnique.edu}, \\
Jean-Michel Lasry \thanks{Universit{\'e} Paris-Dauphine, France},
Pierre-Louis Lions\thanks{Coll{\`e}ge de France,  Paris, France } ,
Antoine Rostand\thanks{Kayrros  Paris, France. a.rostand@kayrros.com},
Jos{\'e} Scheinkman\thanks{Columbia University, New-York, USA. js3317@columbia.edu }
}
\begin{document}
\maketitle

\begin{abstract}
We  propose a plausible mechanism for the short term dynamics
of the oil market based on the interaction of economic agents, namely a major agent (a monopolistic cartel),
a fringe of competitive producers and a crowd of arbitrageurs who store the resource. This model
is linked to a separate work by the same authors on a long term model for the oil industry.
The present model leads to a system of two coupled  non linear partial differential equations, with a new type of  boundary conditions linked to constraints on the storage capacity.  These boundary conditions play a key role and translate the fact that when  storage is either full or empty, the cartel has an enhanced strategic power, and may tune the price of the resource. The model is discussed in details, as well as some simpler variants. A finite difference scheme is proposed and numerical simulations are reported. The latter result in apparently surprising
facts:  1)  the optimal control of the cartel (i.e. its level of production) is  a discontinuous function of the state variables; 2) there is a cycle, which takes place around the shock line. We discuss these phenomena in details, and show that they may explain what happened  in 2015 and 2020.
\end{abstract}

\section{Introduction}\label{sec:introduction}

For several decades, the oil industry has had  structural characteristics that make it exceptional from the viewpoint of economic theory.
A main feature is the existence of a  monopolistic cartel of producers lasting since 1960. This is partly explained by the fact that 
international law does not prohibit  a worldwide coalition, while 
in most countries, the law forbids lasting and non regulated nationwide monopolies.
Regardless of legal issues, the existence of a lasting equilibrium between a cartel
and a competitive fringe of producers may arise only in particular economic situations.
In a separate work in progress, see \cite{edmond1}, we have proposed an economic model  that we have named  {\sl Edmond 1}\footnote{  after the famous writer Edmond Rostand (\cite{cyrano}), a relative of  one of the authors}, that is motivated by the following observations:
\begin{itemize}
\item First, oil reserves and  production capacity \footnote{In our model, the production of the competitive fringe is always the same fraction of the capacity of production, which is of course an approximation of what actually happens in the industry. Thanks to this proportionality, we can use as a state variable the annual production $z$ of the fringe. }  have been driven for several decades by production and investment. Of course, production makes reserves decrease; conversely, new reserves can be found and new production facilities can be built.
  Moreover, for most oil fields, 
  the necessary investments  from  prospecting to   construction of production facilities ({\sl CAPEX}), have been rather stable for several decades, of the order of 10 to 20 dollars per barrel eventually produced . Due to inefficiencies that will be discussed later, the flow of CAPEX, i.e. the CAPEX per year, has a decreasing yield. Indeed, it is observed that the yearly increase in production capacity  is proportional to  
the square root of  CAPEX flow .%
\footnote{This means that if  during a given year, the flow of CAPEX is of $300$ Bn \$  and increases the  production capacity by say $30$ Bn barrels, then, during the same year,  a flow of CAPEX of $600$ Bn \$ would have increased the capacity by $\sqrt{2}\times 30= 1.41\times 30$ Bn barrels. }

  Notice that this first observation leads us to
 a radical departure from Hotelling's framework. Indeed, on the  time scale considered here (i.e.  decades, not centuries), it appears that new reserves are the results of research efforts and investments that are endogenously determined; therefore,  exhaustibility of the resource is not relevant.
  Models for mining industries which do not belong to the Hotelling framework have already been  studied by some of the authors, see \cite{MR3575616}.

\item The second observation that explains the equilibrium between the cartel and the competitive fringe is that  the  investments on capacity  by the fringe producers are limited by credit constraints. Indeed,  most producers in the competitive fringe can invest only from their own cash flows; more precisely, 
  almost all fringe producers can invest only if the price of the barrel is high enough, i.e. above a given threshold, and  their investments increase linearly with respect to the price of the barrel when it is above the threshold.

This constraint on investments is mostly due to a set of major risk factors (operational risks, duration of projects, sovereign risk, etc...). This gives an important strategic lever to the monopolistic cartel, all the more important given that the elasticity of demand is extremely low. Indeed,
  \begin{itemize}
  \item In the short run, the price elasticity is around $0.04 \%$, so that a small increase in the cartel's production drives down prices by a substantial amount. This stops investment by the fringe and drives down their capacity, allowing the monopolist to increase market share
\item Conversely, when the  cartel reduces its production by a relatively small amount,  the price increases substantially and the fringe producers invest. However, since  all fringe producers are investing at the same time,  the CAPEX required to obtain an additional unit of capacity is high. In other words, these simultaneous investments by the fringe producers are inefficient.

  This gives the cartel more strategic power and allows the cartel to make additional profits for a rather long period.
  \end{itemize}
\end{itemize}
The two observations  are key ingredients for the equilibrium described in the long term model \cite{edmond1}, that exhibits   fluctuations between high prices that lead to immediate benefits but also a decrease in the cartel's future market shares and low prices which have the opposite effects. 

The model in  \cite{edmond1} is stationary except for a single multiplicative factor on the demand curve.  A numerical calibration  shows that it  is capable of matching key moments of the observed time series in the oil market

\bigskip

Our long term model \cite{edmond1} leaves aside
detailed interactions between the producers and the business of oil storage: it uses a very simple approximation
to the behavior of the arbitrageurs  who buy, store and then sell 
the resource. This simple approximation  is  sufficiently good  fot studying long term dynamics over decades. Yet, at the time scale of a few years, the interaction of the monopolistic cartel with the storage business leads to particular strategic effects.
The aim of the present paper is to propose a model named {\sl  Edmond 2}, which focuses on these interactions in order to shed some light on unusual but striking  features of the oil markets' prices, production and investments.
 The new model,  otherwise as in \cite{edmond1} produces an equilibrium which in some aspects is quite distinct \footnote{The strategic interactions are quite different from those described in the long term model \cite{edmond1}. Therefore, the two models should be considered as complementary insights on the oil markets. We have not mixed them in order to avoid that a complexity that would hide the different insights.} - in particular it can generate very large and brutal changes,
 namely discontinuities in prices, and in the optimal strategies of the cartel as functions of the state variables, in the absence of noise. 
\\
We will consider a large number of arbitrageurs who behave competitively. As is often the case in economics, a competitive set of arbitrageurs  limits the freedom of the cartel to implement strategies that depend on price changes as in \cite{edmond1}. But, in the oil industry, arbitrageurs have their own limits. Indeed, a key element in the present model  is  the existence of capacity constraints (there is a minimum and a maximum capacity) on storage.
These one-side constraints  allow the cartel more strategic power when the capacity limits  are reached. 
These key features linked to the arbitrageurs' constraints result in three different regimes for the slope of the futures curve, namely {\sl 'contango', 'backwardization' and 'standard'}, as will be seen in Section~\ref{sec:numer-simul}.

\bigskip

The paper is organized as follows: the main part of  the new model  is discussed in Section \ref{sec:syst-part-diff}; mathematically, it leads to a system of partial differential equations. Section \ref{sec:syst-part-diff} also includes the description of simpler variants. The issue of the boundary conditions corresponding to situations when the storage facilities are either empty or full is particularly delicate: it is dealt with in Section \ref{sec:boundary-conditions}.  To the best of our knowledge, such boundary conditions appear to be completely new. Their mathematical analysis seems quite delicate, and we give partial theoretical results in Section \ref{sec:few-results-math}.
Section \ref{sec:appr-finite-diff} is devoted to a finite difference method for solving the system of PDEs supplemented with the above mentioned boundary conditions. Finally, numerical simulations are reported in Section \ref{sec:numer-simul}: in particular, we find  large shocks in the oil price. We then comment on the simulations, and stress the qualitative agreement with what occurred at least twice in the recent history,  2015 and 2020.

\section{The models and the systems of partial differential equations}\label{sec:syst-part-diff}

We consider a cartel producing a natural resource 
and facing both a competitive fringe of small producers and a competitive set of arbitrageurs. 
Even though our motivation is to understand some aspects of the oil market in the short or middle term (of the order of few years)
and we may sometimes use a terminology linked to the oil industry (for example, oil for the resource, OPEC for the cartel),
the models proposed below may be applied in  other settings.
\\ 
There are four types of agents:  consumers, a cartel or major agent (OPEC), minor producers forming a competitive fringe, and the arbitrageurs that buy, store and then sell  the resource. The arbitrageurs will most often limit price changes,
 but when  storage capacity limits are reached, the strategic power of the cartel increases dramatically. When no resource is stored,   the cartel has the power to drive  prices up by cutting  production; conversely, when storage is at its maximal level, the  cartel can drive  prices down by increasing  production.
 \\
 For simplicity, we will not consider  the decision making process of the competitive fringe of small producers but 
 instead assume that the dynamics of their global production rate is given as a function of the current state of the world;
 see the above remarks on credit constraints.
\\
In contrast to the long-run analysis in  \cite{edmond1}, our aim is to model the dynamics in horizons of the order of a few years. For this reason, we will take the demand function by final consumers as constant. This assumption is a good approximation for  oil markets,  since, except for the occurrence of major unexpected shocks - such as the arrival of the pandemic in 2020 - there is little variation in the demand function on a yearly horizon (see paragraph ~\ref{sec:disc-what-happ} for a discussion). Nonetheless we will  comment on  the effect of \emph{rare disasters} on equilibria; see paragraph~\ref{sec:disc-what-happ}.
\\
The class  of models described below involves two state variables: 
\begin{enumerate}
\item  the level $k$ of speculative storage 
\item  the level $z$ of  aggregate output of the competitive fringe. 
\end{enumerate}
While $k$, i.e. the level of speculative storage, takes its value in a given interval, say $[k_{\min},k_{\max}]$,
the second state variable, i.e. $z$, may be either discrete or continuous. The physical constraints on the storage capacity will play a key role. Indeed, it will be  shown that in some situations and when the storage level is either minimal or maximal,  the cartel directly controls the price of the resource. \\
Mathematically, all the variants of the model lead to systems of partial differential equations coupling a Hamilton-Jacobi-Bellman equation for the cartel (major agent) and an equation of the type ``master equation''  for the price of the resource, see \cite{PLL,MR3967062,MR3941633}. Note that in the present case, the master equation does not model a crowd of players as in mean field games, but rather an equilibrium reached by a crowd of arbitrageurs. It seems to be the first example in which the master equation does not involve the value of a game between competitive agents, but rather a price fixed by an arbitrage relationship. Other examples will be supplied in forthcoming papers. To the best of our knowledge, the boundary conditions arising from the state constraints are completely new as well.

\bigskip

Our model has several variants,  mathematically simpler but more limited from the modeling viewpoint:
\begin{enumerate}
\item In the first variant, the global production rate $z_t$ of the competitive 
fringe takes its values in the interval $[z_{\min}, z_{\max}]$
\item In the second variant, $z_t$ can take a finite number of values $z_j$,  $j=0,\dots, J-1$, for 
 a positive  integer $J$ (we will only discuss the cases when $J=2$ and $J=1$).
\end{enumerate}

\subsection{A model with two continuous state variables}
\label{sec:model-with-two}

\subsubsection{The dynamics of $k_t$ and $z_t$}
\label{sec:dynamics-k_t-z_t}
The global production rate  $z_t$ of the competitive fringe is assumed to follow the  dynamics 
\begin{equation}
  \label{eq:dynz}
  dz_t= b(k_t,  p_t) dt,
\end{equation}
where $p_t$ stands for the unitary price  of the resource  and $b: [k_{\min}, k_{\max}]\times \R_+ \to \R$ is a given smooth function, which we take of the form
\begin{equation}
  \label{eq:bbb}
b(k,p) =\kappa (\lambda p - \mu) + f(k).
\end{equation}
The first term in (\ref{eq:bbb}), i.e. $\kappa (\lambda p - \mu)$ expresses the direct impact of prices on investments,
 hence on production capacity, as it has been  explained in the introduction. 
A typical possible choice for the second term  in (\ref{eq:bbb}) is 
\begin{displaymath}
  f(k)=   a_1 \left(\frac {k_{\max} -k} {k_{\max} -k_{\min}}\right)^2   -a_2 \left(\frac {k -k_{\min}} {k_{\max} -k_{\min}}\right)^2 ,
\end{displaymath}
with suitable positive constants $a_1$ and $a_2$, in such a way that $f$ has a significant effect on $b$ only for values of $k$ close to  $k_{\min}$ or $k_{\max}$.  The term $f(k)$ can therefore be seen as a modulation of the first term in $b$ near the limits of $k$.
It is a proxy for the time delays between the investment decisions of the producers belonging to the competitive fringe and the actual creation of new capacities of production.
Indeed,  it would  be a strong simplification to assume that the investment has an {\sl instantaneous} effect. This simplification  proves acceptable in a  very long term model  like \cite{edmond1}, i.e. at the scale of decades; it permits to keep the theoretical complexity of the latter model at a reasonable level.  The situation is much different when one deals with the short term, of the order of a few years: it seems necessary to model inertia effects, memory effects and anticipations of delays  between investments and the creation of capacity.
Nevertheless, since we wish to keep the model as simple as possible, we limit ourselves to a proxy
when addressing the   delay effects. 
The  function $f(k)$ accounts for a little increase (respectively decrease) in  production capacities when the storage facilities are close to empty (respectively full). Indeed, close to empty storage  must follow a period when the price is  high, thus the investments of the producers in the competitive fringe are at a high  level; the latter result in an increase of production capacity, i.e. an increase of $z$, even if the instantaneous price has decreased. The mechanism has to be reversed when the storage facilities are close to full. Hence, we choose $f$ which takes a value  of the order of  $1\%$ for $k\sim k_{\min}$ and  $-1\%$ for $k\sim k_{\max}$.
Note that  an accurate model with more state variables is possible, but  it would be more difficult to understand and to simulate numerically
 (in particular because it would increase the dimensionality of the problem).
\begin{remark}\label{sec:dynamics-k_t-z_t-1}
 We may add some  in the dynamics of $z_t$ and replace \eqref{eq:dynz} with
 \begin{equation}
   \label{eq:92}
     dz_t= b(k_t,  p_t) dt + \sqrt{2\nu_z} dB_t,
 \end{equation}
 where $B_t$ is a Brownian motion. This randomness brings the model closer to reality, since a large number of various shocks with relatively small amplitudes  arise in the whole industry, from production to demand. 
\end{remark}

\begin{remark}\label{sec:dynamics-k_t-z_t-2}
For the numerical simulations, we have to limit ourselves to $z\in [z_{\min}, z_{\max}]$ for  $0< z_{\min}\le z_{\max}$.  The bounds $z_{\min}$ and $z_{\max}$ are chosen from historical data. Hence, it is convenient to add to $b$ another term $\widetilde b(z)$ which vanishes away from $z=z_{\min}$ and $z=z_{\max}$, and which is negative near  $z=z_{\max}$ and positive near  $z=z_{\min}$. The role of this term is purely technical, and it does not affect the numerical results far enough from  $z=z_{\min}$ and $z=z_{\max}$. Therefore, we do not wish to lean too much on it. However, this explains why, in what follows,  $b$ may also depend explicitly on $z$.
\end{remark}

Let us compare the choice of $b$ made in \eqref{eq:bbb} with the one made in the long term model \cite{edmond1}. In \cite{edmond1} (with a very simplified description of the storage business), we take $b$ as follows:
\begin{equation}
  \label{eq:bbb2}
b(z,p) =- a z + \lambda p - \mu  ,
\end{equation}
with $\lambda\approx 0.4$, $\mu\approx 10$ et $a\approx 10\%$. The fact that  the storage $k$ is not involved in (\ref{eq:bbb2}) obviously comes from the fact that  it is not a state variable in \cite{edmond1}.
 On the other hand,  observe that in the oil industry, the range $k_{\max}-k_{\min}$ of the storage available for  arbitrageurs is relatively small (of the order of $5-7\%$ of  annual production, see Section~\ref{sec:numer-simul}). Given this small range and the small elasticity of demand (see the next paragraph), the strategy of the cartel can be implemented by tuning its production $q_t$ within a small range of values. Indeed, it has been  observed for decades that the {\sl spare production capacity} of the cartel mostly varies between $3\%$ and $5\%$
and that the market share $q_t$ of the cartel stays close  to $42\%$. This aspect is taken into account in \cite{edmond1} and explained immediately after equation~(\ref{eq:93}) in  paragraph \ref{sec:equilibrium}. But, to keep the present model focused on the interactions between the cartel and the arbitrageurs, we will describe the range of possible production strategies of the cartel  by a proxy, namely a cost for the deviation of the cartel's production $q_t$ from a target production $q_\circ\approx 42\%$; this cost will have  the form $\alpha(q_t-q_\circ)^2/2$
and be a part of  the running cost in the optimal control problem solved by the cartel.
As a consequence, $z_t$ will stay close to $z_{\circ}\approx 58\%$. Therefore, in the present model,  since the time scale is of the order of a few years, we may neglect the variations of $z$ in (\ref{eq:bbb2}) because they are small on this time scale (these small variations would induce only a small correction in the strategy of the cartel), i.e. replace the term $az$ by $a z_{\circ}$.

\medskip

 The demand of the consumers is a decreasing function of the price of the resource; after a suitable choice of units, the simplest  demand function is 
 \begin{displaymath}
 D(p)=1-\epsilon p,
 \end{displaymath}
 where the parameter $\epsilon$ stands for the elasticity of demand. Note that it would be more appropriate to set
 $D(p)=\max(0,1-\epsilon p)$, but in the regime that will be considered, the price $p$ will never exceed $1/\epsilon$. This linear demand function  fits well the observed data in the usual range of oil prices.
 \\ 
 The control variable of the monopoly is its production rate $q_t$.  Matching demand and supply yields
 $dk_t=  (q_t + z_t -D(p_t) )dt$.
 However, we will instead consider a slightly more general dynamics of $k_t$, possibly including  some   small noise in  storage capacities:
 \begin{equation}
   \label{eq:3}
 dk_t=( q_t + z_t-D(p_t)) dt +  \sigma(k_t) dW_t,  \quad \hbox{for }  k_{\min}\le k_t \le k_{\max},
 \end{equation}
 where $(W_t)$ is a Brownian motion.  We suppose that the volatility $k\mapsto \sigma(k)$ is a smooth nonnegative function that vanishes
  at $k=k_{\min}$ and $k=k_{\max}$ and that the quantities $\frac {\sigma(k)}{k-k_{\min}}$ and $\frac {\sigma(k)}{k_{\max}-k}$ are bounded. This assumption will play an important role in the discussion of the boundary conditions that will be made in paragraph \ref{sec:boundary-conditions} below.

\subsubsection{Equilibrium}\label{sec:equilibrium}
We look for a stationary equilibrium.  Given the unit price of the resource, the cartel solves an optimal control problem. 
 Let $(k,z)\mapsto U(k,z)$ be the associated value  function. The price, described by a function $p(k,z)$, is fixed by ruling out 
 opportunities for arbitrage.  We will see that the functions $U$ and $p$ satisfy a system of two coupled partial differential equations.\\
The optimal control problem solved by the cartel knowing the trajectory of $p_t$ is:
\begin{equation}
  \label{eq:93}
  U(k,z)=\sup_{q_t} \EE\left( \left.\int_0^\infty e^{-rt}    \left( (p_t-c)q_t -\alpha \frac {|q_t-q_\circ|^2}2   \right)   dt\right|  (k_0,z_0)=(k,z) \right),
\end{equation}
where $r$ is a positive discount factor, $c$ is the cost related to the  production of a unit of resource.
To complement what  has already been said in paragraph~\ref{sec:dynamics-k_t-z_t},  the   penalty term
$\alpha \frac {|q_t-q_\circ|^2}2$ is introduced because
 for decades, the interactions between the cartel and the competitive fringe have resulted in the fact
 the production level of the cartel has oscillated around $q_0\approx 42\%$ with a standard deviation of a few $\%$. These interactions have been  modeled  in \cite{edmond1} and briefly described in the introduction. They partly rely on the fact that the cartel can tune its spare production capacities (in a range varying from $3\%$ to $5\%$) in order  to drive  prices up or down and possibly to deter the fringe from investing.
Here, as already explained in paragraph~\ref{sec:dynamics-k_t-z_t},
this aspect is  taken into account via a proxy, namely the contribution $\alpha \frac {|q_t-q_\circ|^2}2$ to the cost.

The dynamic programming  principle yields that the value function is a solution of the following Hamilton-Jacobi-Bellman equation:
\begin{equation}\label{eq:HJB1}
 -r U  + \sup_{q\ge 0}  \left( -\alpha \frac {|q-q_\circ|^2} 2 +(p-c) q + \left(q+z-D\left(p\right)\right) \partial _k U \right)  + b(k,z,p) \partial_z U +\frac {\sigma^2(k) }2  \partial_{kk} U =0.
\end{equation}
Introducing the Hamiltonian
\begin{equation*}
  \begin{split}
  H(z,p, \xi)   &=\sup_{q\ge 0}     \left( -\alpha \frac {|q-q_\circ|^2} 2 +(p-c) q + \xi\left(q+z-D\left(p\right)\right)  \right)
  \end{split}
\end{equation*}
in which the maximum is reached by $q^*= \max(0,q_\circ + \frac 1 \alpha (p-c+\xi) )$, 
(\ref{eq:HJB1}) can be written:
\begin{equation*}
 -r U  + H(z,p, \partial _k U)   + b(k,z,p) \partial_z U +\frac {\sigma^2(k) }2  \partial_{kk} U =0.
\end{equation*}
Since, in the regime that will be considered, $q_\circ + \frac 1 \alpha (p-c+\partial_k U) )  $ will always be nonnegative,
we omit for simplicity the constraint $q\ge 0$ in the definition of the Hamiltonian: hereafter, we set 
\begin{equation}\label{eq:4}
  \begin{split}
  H(z,p, \xi)   &=\sup_{q\in \R}     \left( -\alpha \frac {|q-q_\circ|^2} 2+(p-c) q +\xi \left(q+z-D\left(p\right)\right)  \right)\\ & = 
\frac 1 {2\alpha}  (p-c+\xi)^2  +\xi(z-D(p))+q_\circ (p-c-\xi),  
  \end{split}
\end{equation}
and the optimal production rate at $k_t=k$ and $z_t=z$ is given by the feedback law
\begin{equation}\label{eq:5}
  \begin{split}
q^*(k,z) &=   D_\xi H\left(z,p(k,z),  \partial _k U (z,p)\right)-z +D(p(k,z))\\
&= q_\circ + \frac 1 \alpha (p(k,z)-c+\partial_k U(k,z)) ).
  \end{split}
\end{equation}
Let us now turn to the price of a unit of resource : ruling out opportunities for arbitrage
implies that the price process obeys the following relation,
\begin{displaymath}
  p_t= \EE\left( \left. e^{-r \delta t } p_{t+\delta t} - \int_{t}^{t+\delta t} g(k_s) ds \right | (k_t,z_t)    \right),
\end{displaymath}
where $g(k)$ is the cost of storing  a unit of resource per  unit of time when the level of storage is $k$. Recalling that $p_t=p(k_t,z_t)$,  Ito formula yields: 
\begin{displaymath}
  -r p  +  D_\xi H(z,p, \partial_k U) \partial_k p + b(k,z,p) \partial_z p  +\frac {\sigma^2(k) }2  \partial_{kk} p -g(k)=0.
\end{displaymath}
To summarize,  the system of PDEs satisfied by $U,p$ is 
\begin{eqnarray}
  \label{eq:6}
0&=& -r U + H(z,p, \partial_k U) + b(k,z,p) \partial_z U +\frac {\sigma^2(k) }2  \partial_{kk} U , \\
\label{eq:7}
0&=& -r p  +  D_\xi H(z,p, \partial_k U) \partial_k p + b(k,z,p) \partial_z p +\frac {\sigma^2(k) }2  \partial_{kk} p -g(k),
\end{eqnarray}
for $k_{\min}<k<k_{\max}$ and $z_{\min}<z<z_{\max}$ and with $H$ given by (\ref{eq:4}).
We will see in Section \ref{sec:numer-simul} below  that  (\ref{eq:6}-\ref{eq:7}) may have singular solutions,
which consist of a discontinuous dynamics  and price discontinuities, which are actually observed in the historical data (see paragraph \ref{sec:numer-simul})
\begin{remark}
  The  system   \eqref{eq:6}-\eqref{eq:7} must be completed with boundary conditions. We will see that the  boundary conditions at $k=k_{\min}$ and $k=k_{\max}$ are extremely unusual from the mathematical point of view and have  important economic
  implications.
\end{remark}

Note that  equation (\ref{eq:7}) is nonlinear with respect to $p$. It is reminiscent of the master equations discussed in \cite{MR3941633}.\\
 Note also that it seems  possible to refine the present model by considering that the arbitrageurs running the speculative storage business are rational agents playing a mean field game. This would lead to a more involved model of a mean field game with a major agent, see \cite{MR3851543}. Yet,  the resulting system of partial differential equations would have the same structure as  (\ref{eq:6})-(\ref{eq:7}).

\subsection{A variant in which the production of the fringe is a two-state Poisson process}
\label{sec:prod-fringe-two}


We now  introduce a variant which aims at keeping the essential features of the previous model while being simpler mathematically. 
We consider a situation in which the production rate $z_t$ can take only two values $0\le z_0<z_1$ and is described by  
 a stochastic Poisson process  with intensities that may depend on $k_t$ and $p_t$:
\begin{equation}
  \label{eq:8}
  \left\{\begin{array}[c]{rcl}
    \PP\left(z_{t+\dt} =z_0\right| \left.  z_{t} =z_0 \right)&=& 1- \lambda_0(k_t,p_t) \dt +o(\dt), \\ 
    \PP\left(z_{t+\dt} =z_1\right| \left.  z_{t} =z_0 \right)&=& \lambda_0(k_t,p_t) \dt +o(\dt),  \\
    \PP\left(z_{t+\dt} =z_1\right| \left.  z_{t} =z_1 \right)&=& 1-\lambda_1(k_t,p_t) \dt +o(\dt), \\
    \PP\left(z_{t+\dt} =z_0\right| \left.  z_{t} =z_1 \right)&=& \lambda_1(k_t,p_t) \dt +o(\dt).
  \end{array}\right.
\end{equation}
All the other features of the model are the same as in paragraph \ref{sec:model-with-two}, in particular, 
the dynamics of $k_t$ is still given by (\ref{eq:3}).  The optimal value of the cartel and the price are described by 
  $U(k,z_j)=U_j(k)$ and $p(k,z_j)=p_j(k)$,
 where for $j=0,1$, the  real values functions $U_j, p_j$  are defined on $[k_{\min}, k_{\max}]$
and satisfy a system of four coupled differential equations. 
\\
Introducing the Hamiltonians 
\begin{displaymath}
  H_j(p, \xi)= 
\frac 1 {2\alpha}  (p-c+\xi)^2  +\xi(z_j-D(p))+q_\circ (p-c-\xi),  
\end{displaymath}
(we still omit the constraint that the production rate is nonnegative),
 and repeating the arguments contained in paragraph~\ref{sec:equilibrium}, 
we get the following system of differential equations: 
\begin{eqnarray}
\label{eq:9}
0&=& -r U_j + H_j(p_j, U'_j) +\lambda_j(k,p_j) (U_\ell-U_j) +\frac {\sigma^2(k) }2   U''_j, \\
\label{eq:10}
0&=& -r p_j  +  D_\xi H_j(p_j, U'_j) p_j'  +\lambda_j(k,p_j) (p_\ell-p_j)-g(k)+\frac {\sigma^2(k) }2 p''_j ,
\end{eqnarray}
 for $j=0,1$, $\ell=1-j$, and $k\in (k_{\min},k_{\max})$. The optimal drift of $k_t$ in (\ref{eq:3}) is then given by
 $D_\xi H_j(p_j(k_t), U'_j(k_t))= \frac 1 {\alpha} (U'_j(k_t) -c +p_j(k_t))   +   z_j-D(p_j(k_t)) +q_\circ$  if $z_t= z_j$.
 \begin{remark}
   Here again, the boundary conditions will be important and non standard.
 \end{remark}
\subsection{  An even simpler model}
\label{sec:simpler-model-which}
It is possible to  simplify further the model by assuming that the production rate of the competitive fringe is a constant $z$.
Introducing the Hamiltonian
\begin{displaymath}
  H(p, \xi)=  \frac 1 {2\alpha}  (p-c+\xi)^2  +\xi(z-D(p))+q_\circ (p-c-\xi),  
\end{displaymath}
and repeating the arguments contained in paragraph~\ref{sec:equilibrium}, 
we get the following system of two differential equations: 
\begin{eqnarray}
\label{eq:11}
0&=& -r U + H(p, U') +\frac {\sigma^2(k) }2   U'', \\
\label{eq:12}
0&=& -r p  +  D_\xi H(p, U') p' -g(k)+\frac {\sigma^2(k) }2 p'' ,
\end{eqnarray}
 for  $k\in (k_{\min},k_{\max})$.

\section{Boundary conditions}\label{sec:boundary-conditions}
The systems of partial differential equations must be supplemented by boundary conditions. We are going to discuss the boundary conditions  for the full model, 
 and more briefly for the simplified variants proposed in \ref{sec:prod-fringe-two} and \ref{sec:simpler-model-which}.

These boundary conditions will translate mathematically the main change in strategic power created by the constraints: while the price is determined by the arbitrageurs when  storage  is neither full nor empty, the price is driven by the cartel  when  storage is full or  empty. As always in partial differential equations, these boundary conditions are key for the determination of the solution,
 and therefore of the behaviours of both the cartel and the arbitrageurs. For example, if the range $k_{\max}-k_{\min}$ was very small, then the solution would be mostly determined by the boundary conditions, which translates the fact that  the cartel could neglect the impact of the arbitrageurs. 
\subsection{The boundary conditions associated with the model discussed in paragraph \ref{sec:model-with-two} }
\label{sec:bound-cond-assoc}

\subsubsection{Boundary conditions at  $z=z_{\min}, z_{\max}$}
\label{sec:bound-cond-at}
No boundary conditions are needed at $z=z_{\min}$ and $z=z_{\max}$, because of the assumptions on $b$.

\subsubsection{Boundary conditions at  $k=k_{\min}$ and $k=k_{\max}$}
\label{sec:bound-cond-at-1}
\paragraph{Preliminary : monotone envelopes of $\xi\mapsto H(z,p, \xi)$.}
For describing the boundary conditions linked to the state constraints $ k_{\min}\le k_t\le k_{\max}$, 
it is useful to introduce the nonincreasing and nondecreasing envelopes of the function $\xi\mapsto H(z,p, \xi)$: we set
\begin{equation}
  \label{eq:13}
 H_{\downarrow}(z,p,\xi)= \max_{q\le D(p)-z} \left( -  \frac \alpha 2 (q-q_\circ)^2   + (p-c) q   + \xi (  q+z -D(p) \right),
 \end{equation}
and 
\begin{equation}
    \label{eq:14}
  H_{\uparrow}(z,p,\xi)= \max_{q\ge D(p)-z}  \left( -  \frac \alpha 2 (q-q_\circ)^2   + (p-c) q   + \xi (  q+z -D(p) \right) .
\end{equation}
The Hamiltonian $H_{\downarrow}(z,p,\xi)$  (resp.  $H_{\uparrow}(z,p,\xi)$)  corresponds 
 to the controls $q$ such that the drift of $k_t$ in (\ref{eq:3}) is nonpositive (resp. nonnegative).
It may also be convenient to set 
\begin{equation}
  \label{eq:15}
H_{\min}(z,p)= \min_{\xi} H(z,p,\xi)= -\frac \alpha 2 \left( D(p)-z-q_\circ\right)^2 +(p-c)(D(p)-z)
\end{equation}
which corresponds to the control $q=D(p)-z$ for which  the drift of $k_t$ in (\ref{eq:3}) vanishes.
Note that $p\mapsto H_{\min}(z,p)$ is strongly concave with respect $p$.
 It is easy to check that 
\begin{displaymath}
  H(z,p,\xi)=H_{\downarrow}(z,p,\xi)+H_{\uparrow}(z,p,\xi)-H_{\min}(z,p).
\end{displaymath}
The optimal values of $q$ in the definition of  $H_{\downarrow}(z,p,\xi)$ 
and  $H_{\uparrow}(z,p,\xi)$ are
\begin{eqnarray}
  \label{eq:16}
q^*_{\downarrow}(z,p,\xi)=\min\Bigl(D(p)-z,q_\circ+\frac  {p-c+ \xi} \alpha \Bigr) ,\\
\label{eq:17}
q^*_{\uparrow}(z,p,\xi)=\max\Bigl(D(p)-z,q_\circ+\frac  {p-c+ \xi} \alpha \Bigr) .
\end{eqnarray}
Hence,
\begin{eqnarray}
\label{eq:18}   H_{\downarrow}(z,p,\xi)=
\frac 1 2\left( \left(  \sqrt \alpha    (z   -D(p) +q_0)      +      \frac 1 {\sqrt\alpha}   (p-c+\xi) \right)_- \right)^2  +H_{\min}(z,p),
\\
\label{eq:19}   H_{\uparrow}(z,p,\xi)=
\frac 1 2\left( \left(  \sqrt \alpha    (z   -D(p) +q_0)      +      \frac 1 {\sqrt\alpha}   (p-c+\xi) \right)_+\right)^2 +H_{\min}(z,p).
\end{eqnarray}

\begin{assumption}
  \label{sec:prel-:-monot}
Hereafter, we assume that for all $k\in [k_{\min}, k_{\max}]$, $z\in [z_{\min}, z_{\max}]$, $\xi\in \R$, the function   $p\mapsto H_{\min}(z,p)  + \xi b(k, z,p) $ is strongly concave.
\end{assumption}

\paragraph{Boundary conditions at  $k=k_{\min}$.}

In view of the assumptions made on $\sigma$, it is not restrictive to focus  
on the deterministic case: we take  $\sigma=0$ for simplicity. 
\\
The state constraint $k_t\ge k_{\min}$ implies that $q^*(k_{\min},  z) +z -D(p(k_{\min},  z))\ge 0$. 
Two situations may occur:
\begin{enumerate}
\item If $ \frac {\partial_k U(k,z) -c +p(k,z)} \alpha  +z-D(p)+q_\circ >0 $ for $k$ near $k_{\min}$, then the optimal strategy results in increasing the level of storage. This means that in (\ref{eq:7}), the drift $ D_\xi H(z,p, \partial_k U)$ is positive for $k$  near $k_{\min}$, and  no boundary condition is needed for $p$.
\item On the contrary, if $ \frac {\partial_k U(k,z) -c +p(k,z)} \alpha  +z-D(p)+q_\circ \le 0 $ for $k$ near $k_{\min}$, 
then the optimal drift of $k_t$ in (\ref{eq:3}) must vanish at $k=k_{\min}$, i.e.  $q+z -D(p(k_{\min},z))=0$. 
This relation and the strict monotonicity of $D$ imply that $p$ can  be considered as the control variable at $k=k_{\min}$. 
In other words,  the monopoly directly controls the price in this situation.

On the other hand, to rule out   arbitrage opportunities  
but taking into account the 
 state constraints, it is immediate  
 that the price process satisfies 
\begin{displaymath}
  p_t\ge \EE\left( \left. e^{-r \delta t } p_{t+\delta t} - \int_{t}^{t+\delta t} g(k_s) ds \right | (k_t=k_{min},z_t)    \right).
\end{displaymath}
Since the optimal drift of $k_t$ is $0$, we obtain  
\begin{equation}
  \label{eq:20}
  r p(k_{\min}, z)   - b\Bigl(k_{\min},z, p(k_{\min},z) \Bigr) \partial_z p(k_{\min}, z) +g(k_{\min}) \ge 0.
\end{equation}
Another way to understand (\ref{eq:20}) is as follows:  we expect that, in the present case,  $p$ is nonincreasing with respect to $k$ for $k$ near $k=k_{\min}$. Indeed, if  $p$ was increasing with respect to $k$ for $k$ near $k=k_{\min}$,  then the arbitrageurs would increase the level of storage, i.e. $dk_t$ would be positive, in contradiction with the assumption. Then, plugging this information in (\ref{eq:7}) implies (\ref{eq:20}).

Turning back to the cartel, we deduce from the considerations above that, among the strategies consisting of keeping $k_t$ fixed at $k_{\min}$ for $z_t=z$, the optimal one is 
 \begin{eqnarray}
   q^* &=&  D(p^*) -z, \\
   p^*&=&\argmax_{  \pi:    r\pi\ge    b(k_{\min}, z, \pi) \partial_z p-g(k_{\min})} F(\pi, \partial_z U),
 \end{eqnarray}
 where
 \begin{equation}
   \label{eq:21}
F(\pi, \partial_z U)= H_{\min}(z,\pi)  + b(k_{\min}, z,\pi) \partial_z U,
 \end{equation}
and $H_{\min}(z,\pi)$ is defined in (\ref{eq:15}).
Note that $\pi^*$ is unique from Assumption~\ref{sec:prel-:-monot} and depends on $z,\partial_z U, \partial_z p$.
In this situation, the  nonlinear boundary  condition
\begin{equation}
  \label{eq:22}
  p=  p^*(z, \partial_z U, \partial_z p)
\end{equation}
must be imposed at $(k_{\min}, z)$.
\end{enumerate}

\paragraph{Summary.}
 Setting $\ds p(k_{\min,+},z)=\mathop{\lim}_{k-k_{\min}\to 0+} p(k,z)$,  
another way of formulating the boundary conditions at $k=k_{\min}$ is:
\begin{itemize}
\item The nonlinear condition (\ref{eq:22})
, i.e. 
  \begin{equation}\label{eq:23}
p  =  p^*(z, \partial_z U, \partial_z p),
  \end{equation} 
understood in a weak sense, (i.e. it holds only if the optimal drift  $ \partial_k U(k,z) -c +p(k,z)$
 is $\le 0$ near $k=k_{\min}$), and where $p^*(z, \partial_z U, \partial_z p)$  achieves the maximum in (\ref{eq:26}) below

\item the equation for $U$ can be written
\begin{equation}
\label{eq:24}
-rU+ \max(A, B)=0   ,
\end{equation}
with 
\begin{eqnarray}
    \label{eq:25}
A&=& \ds  H_{\uparrow}(z,p(k_{\min,+},z),\partial_k U)
 + b\Bigl(k_{\min},z,p(k_{\min,+},z)\Bigr) \partial_z   U   ,  \\
  \label{eq:26}
B&=&\ds 
\max_{  \pi:    r\pi\ge    b(k_{\min}, z,\pi) \partial_z p-g(k_{\min})} F(\pi, \partial_z U),
\end{eqnarray}
and $F$ is given by (\ref{eq:21}).
\end{itemize}
Note that, to the best of our knowledge, this set of boundary conditions, associated to the system (\ref{eq:6}-\ref{eq:7}) and to the state constraint $k\ge k_{\min}$, has never been proposed and a fortiori analyzed.

\paragraph{Boundary conditions at $k=k_{\max}$.}
Arguing as above and setting $p(k_{\max,-},z)=\lim_{k-k_{\max}\to 0-} p(k,z)$,
 the boundary conditions at $k=k_{\max}$ can be written as follows:

\begin{itemize}
\item A nonlinear condition for $p$ of the form  
  \begin{equation}
\label{eq:27}
 p=  p^{**} (z,\partial_z U, \partial_z p), \end{equation}
 understood in a weak sense (i.e. it holds only if the optimal drift  $ \partial_k U(k,z) -c +p(k,z)$ is $\ge 0$ near $k=k_{\max}$),
where $p^{**} (z,\partial_z U, \partial_z p)$  achieves the maximum in (\ref{eq:30}) below (it is  unique from Assumption~\ref{sec:prel-:-monot}).
\item An equation for $U$: 
\begin{equation}
\label{eq:28}
-rU+ \max(C, D) =0   ,
\end{equation}
with 
\begin{eqnarray}
\label{eq:29}
C&=& \ds H_{\downarrow}(z,p(k_{\max,+},z),\partial_k U)  + b(k_{\max},z,p) \partial_z   U ,
\\
\label{eq:30}
D&=&\ds \max_{   \pi:    r\pi\le    b(k_{\max}, z,\pi) \partial_z p-g(k_{\max})} G(\pi, \partial_z U),  
\end{eqnarray}
and 
\begin{equation}
  \label{eq:31}
G(\pi, \partial_z U)= H_{\min}(z,\pi)  + b(k_{\max}, z,\pi) \partial_z U.  
\end{equation}
\end{itemize}

\subsection{The boundary conditions associated with the model  in ~\ref{sec:prod-fringe-two}}
\label{sec:bound-cond-assoc-1}
The boundary conditions associated with the system (\ref{eq:9}-\ref{eq:10}) are obtained in the same manner as in the previous case.
To avoid repetitions, we focus on the boundary $k=k_{\min}$, because 
the needed modifications  with respect to paragraph \ref{sec:bound-cond-at-1} are similar for $k=k_{\max}$ and $k=k_{\min}$. 
The interested reader will easily find the boundary conditions at $k=k_{\max}$ from paragraph \ref{sec:bound-cond-at-1} and what follows. \\
As above, we set 
\begin{eqnarray}
\label{eq:32}
H_{j,\min}(p)&=& \min_{\xi} H_j(p,\xi)= -\frac \alpha 2 \left( D(p)-z_j-q_\circ\right)^2 +(p-c)(D(p)-z_j),\\
\label{eq:33}
  H_{j,\downarrow}(p,\xi)&=
& 
\max_{q\le D(p)-z_j} \left( -  \frac \alpha 2 (q-q_\circ)^2   + (p-c) q   + \xi (  q+z_j -D(p) \right)
\\
 \notag  &=& 
\frac 1 2\left( \left(  \sqrt \alpha    (z_j   -D(p) +q_\circ)      +      \frac 1 {\sqrt\alpha}   (p-c+\xi) \right)_- \right)^2  +H_{j,\min}(p),
\\  \label{eq:34}
  H_{j,\uparrow}(p,\xi)&=& \max_{q\ge 1-z_j-\epsilon p}
\left( -  \frac \alpha 2 (q-q_\circ)^2   + (p-c) q   + \xi (  q+z_j -D(p) \right)
\\
\notag  &=& 
\frac 1 2\left( \left(  \sqrt \alpha    (z_j   -D(p) +q_\circ)      +      \frac 1 {\sqrt\alpha}   (p-c+\xi) \right)_+\right)^2 +H_{j,\min}(p).
\end{eqnarray}
The optimal values of $q$ in the definition of  $H_{j,\downarrow}(p,\xi)$ 
and  $H_{j,\uparrow}(p,\xi)$ are
\begin{eqnarray}
\label{eq:35}
q^*_{j,\downarrow}(p,\xi)=\min\Bigl(D(p)-z_j,q_\circ+\frac  {p-c+ \xi} \alpha \Bigr)   ,\\
\label{eq:36}
q^*_{j,\uparrow}(p,\xi)=\max\Bigl(D(p)-z_j,q_\circ+\frac  {p-c+ \xi} \alpha \Bigr) .
\end{eqnarray}

\paragraph{Boundary conditions at $k=k_{\min}$.}
 Setting $\ds p_{\ell,+}=\mathop{\lim}_{k-k_{\min}\to 0+} p_\ell(k)$, $\ell=0,1$, the boundary conditions at $k=k_{\min}$ are as follows: for $i=0,1$ and $j=1-i$,
\begin{itemize}
\item a condition of the form 
  \begin{equation}
    \label{eq:37}
p_j=  p_j^*(U_j, U_\ell, p_{\ell,+}), \quad \hbox{with} \quad \ell=1-j,
  \end{equation}
 understood in a weak sense, (i.e. it holds only if $  U'_j(k) -c +p_j(k)\le 0$ for $k$ near $k_{\min}$),
  where $p_j^{*}(U_j, U_\ell, p_{\ell,+})$  achieves the maximum in (\ref{eq:40}) below (it is supposed to be unique).
\item the equation for $U_j$ can be written
\begin{equation}
\label{eq:38}
-rU_j+ \max(A, B)=0   ,
\end{equation}
with
\begin{eqnarray}
\label{eq:39}
A&=& \ds  H_{j,\uparrow}(p_{j,+},U_j') + \lambda_j(k_{\min},p_{j,+}) (U_\ell-U_j)   ,  \\
\label{eq:40}
B&=&\ds  \mathop{\max}_{  p:  (r+\lambda_j(k_{\min},p)) p - \lambda_j(k_{\min},p) p_{\ell,+} +g(k_{\min}
)\ge 0 } F_j(p, U_j, U_\ell),  
\end{eqnarray}
with $\ell=1-j$, and 
\begin{equation}\label{eq:41}
F_j(p, U_j, U_\ell)= H_{j,\min}(p)  +\lambda_j(k_{\min},p) (U_\ell-U_j).  
\end{equation}
\end{itemize}


\subsection{The boundary conditions associated with the model discussed in paragraph~\ref{sec:simpler-model-which}}
\label{sec:bound-cond-assoc-2}
Here also, we focus on $k=k_{\min}$ to avoid repetitions.
Let us set
\begin{eqnarray}
\label{eq:47}
H_{\min}(p)&=& \min_{\xi} H(p,\xi)=-\frac \alpha 2 \left( D(p)-z-q_\circ\right)^2 +(p-c)(D(p)-z),\\
\label{eq:48}
  H_{\downarrow}(p,\xi)&=&\max_{q\le D(p)-z} \left( -  \frac \alpha 2 (q-q_\circ)^2   + (p-c) q   + \xi (  q+z -D(p) \right)\\
 \notag  &=& 
\frac 1 2\left( \left(  \sqrt \alpha    (z  -D(p) +q_\circ)      +      \frac 1 {\sqrt\alpha}   (p-c+\xi) \right)_- \right)^2  +H_{\min}(p),
\\ \label{eq:49}
  H_{\uparrow}(p,\xi)&=&\max_{q\ge 1-z-\epsilon p}
\left( -  \frac \alpha 2 (q-q_\circ)^2   + (p-c) q   + \xi (  q+z -D(p) \right)
\\
\notag  &=& 
\frac 1 2\left( \left(  \sqrt \alpha    (z  -D(p) +q_\circ)      +      \frac 1 {\sqrt\alpha}   (p-c+\xi) \right)_+\right)^2 +H_{\min}(p).
\end{eqnarray}

\paragraph{Boundary conditions at $k=k_{\min}$.}
 Setting $\ds p_{+}=\mathop{\lim}_{k-k_{\min}\to 0+} p(k)$, the boundary conditions at $k=k_{\min}$ are as follows: 
\begin{itemize}
\item a condition of the form 
  \begin{equation}
\label{eq:50}
p=  p^*,
  \end{equation}
 understood in a weak sense, (i.e. it holds only if $  U'(k) -c +p(k)\le 0$ for $k$ near $k_{\min}$),
  where $p^{*}$  achieves the maximum in  (\ref{eq:53}) below ($p^{*}$  is unique).
\item the equation for $U$ can be written
\begin{equation}
\label{eq:51}
-rU+ \max(A, B)=0   ,
\end{equation}
with 
\begin{eqnarray}
\label{eq:52}
A&=& \ds  H_{\uparrow}(p_{+},U') ,
\\
\label{eq:53}
B&=&\ds  \mathop{\max}_{  p:  r p  +g(k_{\min})\ge 0 } H_{\min}(p) .
\end{eqnarray}
\end{itemize}

\section{Mathematical analysis of the boundary conditions in the one dimensional model}
\label{sec:few-results-math}
In what follows, we discuss how the boundary conditions from   paragraph~\ref{sec:bound-cond-assoc-2} determine uniquely   solution of (\ref{eq:11})-(\ref{eq:12}) near the boundary. Although  the argument proposed below is rather formal,
it  gives useful information on the solutions. More precisely, we are going to see that the boundary conditions induce a  
 unique  expansion of the function $p$ and of the derivative of the value function  $V = \partial_k U$ near the boundary. The system of PDEs satisfied  by $p$ and $V$ is as follows:
\begin{eqnarray}
\label{eq:90}
0&=& -r V + \Bigl(\sigma(k) \sigma'(k)+D_\xi H(p, V)\Bigr) V' +    D_p H(p, U') p'   +\frac {\sigma^2(k) }2   V'', \\
\label{eq:91}
0&=& -r p  +  D_\xi H(p, V) p' -g(k)+\frac {\sigma^2(k) }2 p'' ,
\end{eqnarray}
 for  $k\in (k_{\min},k_{\max})$.
 We focus on the boundary conditions at $k= k_{\min}$. \\
First, in the case in which the drift $D_{\xi}H(p,V)$ is positive near $k_{\min}$,   $p$ and $V$ are expected to be smooth at the boundary.
\\
 Hence, we focus on the case in which the drift $D_{\xi}H(p,V)$ points toward the boundary (i.e. $D_{\xi}H(p,V)\leq 0$).
 We make the following ansatz:
\begin{eqnarray}
  \label{devV}
V(k) &= V(k_{\min}) + \gamma (k - k_{\min})^n + o((k - k_{\min})^n),\\
\label{devp}
p(k) &= p(k_{\min}) - \beta (k - k_{\min})^m + o((k - k_{\min})^m),
\end{eqnarray}
with $n,m \leq 1$. For shortening the notation, let us define the pair $(V_0,p_0) := (V(k_{\min}), p(k_{\min}))$.
\subsection{A singularity is expected}\label{sec:singularity-expected}
Let us explain  why a singular behavior should be  expected near the boundary $k=k_{\min}$. Indeed, assume  that this is not the case 
and that $m=n= 1$; in this situation, from the assumption made on the sign of the drift near the boundary and the constraint $k_t \geq k_{\min}$, we deduce that 
\begin{equation}\label{drift0}
0 = D_{\xi}H(p_0,V_0) =  (\frac{1}{\alpha} + \epsilon) p_0 + \frac{1}{\alpha} V_0 - \frac{c}{\alpha} + z - 1 - q_{\circ}.
\end{equation}
Then, plugging  the ansatz for $V$ and $p$ into (\ref{eq:90})-(\ref{eq:91}) and focusing on the zeroth order terms, we obtain that
\begin{equation}\label{eq:89}
\begin{cases}
rV_0 = \beta(\frac{1}{\alpha}(p_0 - c + V_0) + \epsilon V_0 + q_{\circ}),\\
r p_0 = - g(k_{\min}).
\end{cases}
\end{equation}
The equations in (\ref{eq:89}) and (\ref{drift0}) form a linear system which is over-determined except for a single value of $\beta$.
 Thus, the values of $V_0$, $p_0$ and $\beta$ are determined.  Passing to the first order terms in the expansion of the system, 
we obtain two second order polynomial equations in $\gamma$ and $\beta$, while $\beta$ is already known. It is then easy to observe that
 for a generic choice of the parameters, this system of second order equations is not consistent with the  already obtained 
values of $V_0$, $p_0$ and $\beta$. 
\begin{remark}\label{sec:singularity-expected-1}
Recall that in the case in which the drift is positive near $k_{min}$, no singularity is expected.
\end{remark}

\subsection{Characterization of the singularity}
\begin{proposition}\label{sec:char-sing-1}
If $V$ and $p$ satisfy (\ref{devV}) and (\ref{devp}), then $n = m = 1/2$ and the pair $(V_0,p_0)$ is completely determined
 by the values of $z,\epsilon,\alpha,q_{\circ}$. Moreover, if $(\alpha \epsilon)^2 + \alpha \epsilon > 1$, then there is at most one pair
 $(\gamma,\beta)$ such that (\ref{devV}) and (\ref{devp}) hold.
\end{proposition}
\begin{remark}
\label{sec:char-sing-2}
The latter condition on $\alpha \epsilon$ will be fulfilled in the numerical simulations in \S~\ref{sec:numer-simul} below.  
\end{remark}
\begin{remark}\label{sec:char-sing}
The value of $p_0$ is obviously $p^*$ which has been defined in paragraph~\ref{sec:bound-cond-assoc-2}.
\end{remark}
\begin{proof}Plugging the ansatz into  (\ref{eq:90})-(\ref{eq:91}),
 and using both the boundedness of $g$ and the fact that  $\sigma$ vanishes near $k_{\min}$, 
we deduce that 
\begin{equation}\label{eq:78}
(1 + \alpha\epsilon)V_0 + p_0 = c - \alpha q_{\circ}
\end{equation}
 by identifying the higher order terms in the expansion.
From the state constraint and the sign assumption on the drift, the following also holds:
\begin{equation}\label{eq:79}
(1 +\alpha \epsilon) p_0 +  V_0 = c+ \alpha(1-z ) - \alpha q_{\circ}.
\end{equation}
Since $\alpha \epsilon \notin \{-2; 0\}$, we deduce that 
\begin{equation}\label{eq:80}
\begin{cases}
\ds V_0 = \frac{ z - 1 + \epsilon (c -  \alpha q_\circ) }{\epsilon(2 + \alpha \epsilon)},\\
\ds p_0 = \frac{\epsilon (c- \alpha q_{\circ}) + (1 + \alpha \epsilon)(1 - z)}{\epsilon(2 + \alpha \epsilon)}.
\end{cases}
\end{equation}
Identifying the  higher order terms in the expansion, we see that if $n\ne m$, then $\beta =\gamma = 0$.
 Therefore, $n= m$. Now, if $2n - 1 \notin \{0; 1\}$, identifying the terms of order $2n -1$ leads to 
\begin{equation}\label{eq:81}
\begin{cases}
-(1 + \alpha \epsilon)\beta + \gamma = 0,\\
(1 + \alpha \epsilon)\gamma - \beta = 0.
\end{cases}
\end{equation}
The latter system yields that $\gamma = \beta = 0$. Thus $n \in \{1/2; 1\}$. The only possible value of $n$ is $1/2$ since the case $n = 1$ has 
  already been ruled out. \\
Considering the zeroth order terms, we conclude that
\begin{equation}\label{order0}
\begin{cases}
0 = -rV_0 + \frac{1}{2\alpha}(\gamma - \beta)^2- \epsilon \gamma \beta,

\\
0 = - r p_0 - \frac{\beta}{2\alpha}(\gamma - \beta) + \frac{\epsilon}{2}  \beta^2 - g(k_{\min}).
\end{cases}
\end{equation}
Let us introduce the parameter 
\begin{equation}\label{eq:82}
\lambda = \frac{rV_0}{rp_0 + g(k_{\min})},
\end{equation}
which is well defined since $g\geq 0, p_0 > 0$. Observe that $0\geq \lambda \geq -1$. We deduce that
\begin{equation}\label{eq:83}
\frac{1}{2\alpha}(\gamma - \beta)^2-\epsilon \gamma \beta = - \lambda\frac{\beta}{2\alpha}(\gamma - \beta) + \lambda\frac{\epsilon}{2}  \beta^2,
\end{equation}
then that
\begin{equation}\label{eq:84}
\gamma^2 + (1 - \lambda(1 + \alpha \epsilon))\beta^2 - (1 + 2\alpha \epsilon - \lambda)\gamma \beta= 0 .
\end{equation}
Defining the numbers $x_{\pm}$ by 
\begin{equation}\label{eq:85}
x_{\pm} = \frac{1 + 2\alpha \epsilon - \lambda \pm \sqrt{(1 + 2\alpha\epsilon - \lambda)^2 - 4 (1 - \lambda(1 + \alpha \epsilon))}}{2},
\end{equation}
we finally obtain that 
\begin{equation}\label{eq:86}
(\gamma - x_+ \beta)(\gamma - x_- \beta) = 0.
\end{equation}
Rewriting the second equation in (\ref{order0}), we obtain
\begin{equation}\label{eq:87}
2\alpha g(k_{\min}) + 2\alpha r p_0 = \beta((\alpha\epsilon + 1)\beta - \gamma).
\end{equation}
Thus, $\gamma = x_+ \beta$ is impossible if $x_+ > 1 + \alpha \epsilon$. 
An easy calculation leads to the fact that this last condition is satisfied (i.e. $x_+$ is large enough) if 
\begin{equation}\label{eq:88}
(\alpha \epsilon)^2 + \alpha \epsilon - 1 > 0.
\end{equation}
\end{proof}

\section{Approximation by a finite difference method}
\label{sec:appr-finite-diff}
As already mentioned, the solutions of the system of PDEs may be discontinuous. The numerical scheme must be designed in order to handle these discontinuities. 

Let us focus on the case when 
\begin{equation}
  \label{eq:58}
  \begin{split}
b(k,z,p)\frac {\partial p} {\partial z}&=  (   \phi(k) +\kappa(\lambda p-\mu) )\frac {\partial p} {\partial z}
\\ &=   \frac \partial {\partial z} \left(   \phi(k)  p  +     \frac {\kappa}{2\lambda}   (\lambda p-\mu)^2  \right)  .
  \end{split}
\end{equation}
We are going to use the latter  conservative form in the numerical scheme for (\ref{eq:7}).
Note that 
\begin{displaymath}
     \phi(k)  p  +     \frac {\kappa}{2\lambda}   (\lambda p-\mu)^2  =
 \frac {\kappa \lambda}{2}  \left(p-\frac \mu \lambda +\frac {\phi(k) } {\kappa \lambda}  \right)^2 
-  \frac {\phi^2 (k) } {2 \kappa \lambda}+\frac \mu \lambda  \phi(k)
\end{displaymath}
It is useful to introduce the following numerical flux function:
\begin{equation}
  \label{eq:59}
  \Psi(k, p_\ell,p_r)=\left\{ 
    \begin{array}[c]{ll}
   \ds   \max_{p_\ell\le p\le p_r}\left(   \phi(k)  p  +     \frac {\kappa}{2\lambda}   (\lambda p-\mu)^2  \right),\quad & \ds \hbox{ if } p_\ell\le p_r,\\  
     \ds \min_{p_r\le p\le p_\ell}\left(   \phi(k)  p  +     \frac {\kappa}{2\lambda}   (\lambda p-\mu)^2  \right),\quad &\ds \hbox{ if } p_\ell\ge p_r,\\  
    \end{array}
\right.
\end{equation}
and straightforward calculus leads to 
\begin{equation}\label{eq:60}
      \Psi(k, p_\ell,p_r) =  - \frac {\phi^2 (k)} {2\kappa \lambda} +  \frac \mu \lambda   \phi(k)
 + \frac {\kappa \lambda}{2} \max\left(  
  \left(p_r-\frac \mu \lambda +\frac {\phi(k) } {\kappa \lambda}  \right)_ +^2
,   \left(p_\ell-\frac \mu \lambda +\frac {\phi(k) } {\kappa \lambda}  \right)_ -^2
   \right).      
\end{equation}

Consider a uniform  grid on the rectangle  $[k_{\min}, k_{\max}]\times [z_{\min}, z_{\max}]$: we set  $k_{i}= k_{\min}+i \dk$, $i=0,\dots, N$, with $\dk=\frac {k_{\max}-k_{\min}} N$ and  $z_{j}= z_{\min}+j \dz$, $j=0,\dots, M$, with $\dz=\frac {k_{\max}-k_{\min}} M$
The discrete approximation of $U(k_i, z_j)$ and $p(k_i, z_j)$ are respectively named $U_{i,j}$ and $p_{i,j}$.
\subsection{The discrete version of the system (\ref{eq:6}-\ref{eq:7})}\label{sec:discr-vers-syst}
We use the following notation for the three nodes centered finite difference approximation of 
the second order derivative with respect to $k$:
\begin{displaymath}
  (D^2_{k} U)_{i,j}=   \frac { U_{i+1,j} -2U_{i,j}+U_{i-1,j}}{\dk^2}.
\end{displaymath}
Consider also the first order  one sided finite difference approximations of $\partial_k U$ and $\partial_z U$, namely
\begin{displaymath}
  \begin{split}
  (D_{k, \ell} U)_{i,j}=\frac {U_{i,j}-U_{i-1,j}}{\dk},\quad\quad  &(D_{k, r} U)_{i,j}=\frac {U_{i+1,j}-U_{i,j}}{\dk},\\    
(D_{z, \ell} U)_{i,j}= \frac {U_{i,j}-U_{i,j-1}}{\dz}, \quad \quad  &(D_{z, r} U)_{i,j}=  \frac {U_{i,j+1}-U_{i,j}}{\dz}.
  \end{split}
\end{displaymath}
The advection term with respect to $z$ in (\ref{eq:6}) will be discretized with a first order upwind scheme.
The discrete version of the Hamiltonian $H$ involves the  function $\cH:  \R^4\to \R$,
\begin{displaymath}
  \cH(z,p,\xi_{\ell},\xi_r)= H_{\downarrow}(z,p,\xi_{\ell})+ H_{\uparrow}(z,p,\xi_r)-H_{\min}(z,p),
\end{displaymath}
where $H_{\downarrow}$, $H_{\uparrow}$ and $H_{\min}$ are respectively defined in (\ref{eq:18}), (\ref{eq:19}) and (\ref{eq:15}).
Note that $\cH$ is nonincreasing with respect to $\xi_{\ell}$ and nondecreasing with respect to $\xi_r$.
\\
The discrete version of (\ref{eq:6}) (monotone and first order scheme) is as follows:
 \begin{equation}
\label{eq:61}
0=\left\{
  \begin{array}[c]{l}
\ds -r U_{i,j}  +\frac {\sigma^2(k_i) }2    (D^2_{k} U)_{i,j}\\
\ds +\cH\Bigl(z_j, p_{i, j},   (D_{k, \ell} U)_{i,j}, (D_{k, r} U)_{i,j}  \Bigr) 
\\ \ds  +\max\left(0,b(k_i,z_j,p_{i,j})\right)   (D_{z, r} U)_{i,j} +\min\left(0,b(k_i,z_j,p_{i,j})\right)  (D_{z, \ell} U)_{i,j}
  \end{array}\right. 
 \end{equation}
 for $i=1,\dots, N-1$ and $j=0,\dots, M$. Note that the scheme is actually well defined at $j=0$ (with a slight abuse of notation), because,  since $b(k_i,z_{\min},p_{i,j})\ge 0$, it does not involve $U_{i,-1}$. A similar remark can be made in the case when  $j=M$.
\\
We choose the following discrete version of (\ref{eq:7}):
\begin{eqnarray}
\notag
 g(k_i)&=& -rp_{i,j}+\frac {\sigma^2(k_i) }2    (D^2_{k} p)_{i,j}\\
\label{eq:62}
&&
\ds +
\frac {\partial H_{\downarrow} }{\partial \xi} \Bigl(z_j, p_{i, j},   (D_{k, \ell} U)_{i,j} \Bigr)  (D_{k, \ell} p)_{i,j} 
 +\frac {\partial H_{\uparrow} }{\partial \xi} \Bigl(z_j, p_{i, j},   (D_{k, r} U)_{i,j}  \Bigr)  (D_{k, r} p)_{i,j}
\\ \notag
&& + \frac 1 {\dz} \left( \Psi(k_i, p_{i,j}, p_{i,j+1})-\Psi(k_i, p_{i,j-1}, p_{i,j})\right),
\end{eqnarray}
for $i=1,\dots, N-1$ and  $j=0,\dots, M$. \\
The choice of a monotone scheme for  (\ref{eq:6})-(\ref{eq:7}) will allow us to capture the shocks that have been already mentioned.

 \begin{remark}
  \label{sec:discr-vers-syst-1}
Note that  in  (\ref{eq:6})-(\ref{eq:7}), neglecting all the viscous effects and taking the derivative of  (\ref{eq:6}) with respect to $k$ leads to a  weakly hyperbolic system, (only weakly because the Jacobian matrix has repeated eigenvalues and an incomplete set of eigenvectors, i.e. it is not  diagonalizable). The scheme proposed above may be modified in order to handle what physicists and mathematicians call sonic points and rarefaction waves, by following the ideas in \cite{MR3494346}; in the application presently discussed, it turns out that such an improvement is not necessary
because we did not observe any rarefaction wave.  In a forthcoming paer, we will consider an economic model described by a  weakly hyperbolic system  leading to sonic points, and we will discuss a  numerical scheme  which copes with rarefaction waves.
\end{remark}
\subsection{The discrete scheme at $i=0$}\label{sec:discrete-scheme-at}
In order to write the discrete version of the boundary conditions at $k=k_{\min}$, we set 
\begin{eqnarray}
  \label{eq:64}
A_j&=&\left\{
  \begin{array}[c]{l}
\ds   H_{\uparrow}(z_j,p_{0,j}, (D_{k, r} U)_{0,j} ) \\ \ds
+\max\left(0,b(k_{\min},z_j,p_{0,j})\right)   (D_{z, r} U)_{0,j} +\min\left(0,b(k_{\min},z_j,p_{0,j})\right)  (D_{z, \ell} U)_{0,j}    ,  \end{array}\right.\\
\label{eq:65}
B_j&=&\max_{ p:  K_j( p_{0,j-1},p, p_{0,j+1}) \le g(k_{\min} )} F_{j}(p,U),
\end{eqnarray}
where
\begin{equation}\label{eq:66}
   K_j( p_{0,j-1},p, p_{0,j+1})= -rp + \frac 1 {\dz} \Bigl( \Psi(k_{\min}, p, p_{0,j+1})-\Psi(k_{\min}, p_{0,j-1}, p)\Bigr),
\end{equation}
and 
\begin{equation}
  \label{eq:67}
F_{j}(p,U)= H_{\min}(z_j, p)+  \max\left(0,b(k_{\min},z_j,p)\right)(D_{z, r} U)_{0,j}  +  \min\left(0,b(k_{\min},z_j,p)\right)(D_{z, \ell} U)_{0,j}  .
\end{equation}

The numerical scheme corresponding to the boundary condition  at $i=0$ consists of two equations for each $0\le j\le M$:
\begin{enumerate}
\item  the first equation is \begin{equation}
\label{eq:68}
-rU_{0,j}+ \max(A_j, B_j) =0,
\end{equation}
\item the second equation is either (\ref{eq:69})  or (\ref{eq:70}) below:
 \begin{enumerate}
 \item  if the maximum in (\ref{eq:68}) is achieved by $A_j$, then
 \begin{equation}\label{eq:69}
 g(k_{\min})=\left\{ \begin{array}[c]{l}
\ds 
 -rp_{0,j}  +\frac {\partial H_{\uparrow} }{\partial \xi} \Bigl(z_j, p_{0, j},  (D_{k, r} U)_{0,j}  \Bigr)  (D_{k, r} p)_{0,j}
\\ \ds + \frac 1 {\dz} \Bigl( \Psi(k_{\min}, p_{0,j}, p_{0,j+1})-\Psi(k_{\min}, p_{0,j-1}, p_{0,j})\Bigr)
  \end{array}\right.
\end{equation}
\item otherwise (the maximum in (\ref{eq:68}) is achieved by $B_j$), 
 \begin{equation}
\label{eq:70}
   p_{0,j}= p_j^{*}(U,P),
 \end{equation}
where $ p_j^{*}(U,P)$ achieves the maximum in (\ref{eq:65}).
\end{enumerate}
\end{enumerate}

\begin{remark} There is a unique solution to
  \begin{equation}
\label{eq:71}
 rp=  \frac 1 {\dz} \left( \Psi(k_{\min}, p, p_{0,j+1})-\Psi(k_{\min}, p_{0,j-1}, p)\right) -g(k_{\min}).
\end{equation}
Indeed, (\ref{eq:71}) can be written $\chi(q)= -g(k_{\min})$, with 
  \begin{equation}
\label{eq:72}
  \begin{split}
\chi(q) =&r   \left( \frac {\mu} \lambda -\frac{\phi(k_{\min})}{\kappa \lambda} \right) +rq \\
&-  
 \frac {\kappa \lambda } {2 \dz}
 \max\left( \left(q_{0,j+1}  \right)_+^2 ,  \left(q \right)_-^2   \right)
 +  \frac {\kappa \lambda } {2 \dz} \max\left( (q)_+^2 ,  \left(q_{0,j-1}  \right)_-^2   \right).
  \end{split}
\end{equation}
 The function $\chi$  is increasing  and $\lim_{q\to \pm +\infty}\chi(q)= \pm \infty$. Note also that 
 \begin{displaymath}
   \chi'(q)= r +  \frac {\kappa \lambda } {\dz} q
 \left ( -1_{q\le  -\left(q_{0,j+1}  \right)_+ } +1_{q\ge  \left(q_{0,j-1}  \right)_- }\right) .
 \end{displaymath}
Setting
\begin{displaymath}
  Q=-   \left( \frac {\mu} \lambda -\frac{\phi(k_{\min})}{\kappa \lambda} \right) +   \frac {\kappa \lambda } {2 r\dz}
  \left(q_{0,j+1}  \right)_+^2 
  -  \frac {\kappa \lambda } {2 r \dz}  \left(q_{0,j-1}  \right)_-^2 -\frac {g(k_{\min})}r ,
\end{displaymath}
we see that 
\begin{equation}
  q=\left\{
    \begin{array}[c]{ll}
           \ds  \frac{r -\sqrt{ r^2 +\frac {2 \kappa \lambda}{\dz} \left( r\left(\frac \mu \lambda - 
\frac {\phi(k_{\min})}{\kappa \lambda}\right) +\frac {\kappa\lambda } {2 \dz}  (q_{0,j-1})_- ^2  + g(k_{\min})\right)  }}
{\frac {\kappa\lambda}{\dz} }                 & \ds \hbox{if }\quad Q<-\left(q_{0,j+1}  \right)_+,\\
      Q,\quad & \ds \hbox{if }\quad Q\in \left[ -\left(q_{0,j+1}  \right)_+, \left(q_{0,j-1}  \right)_- \right],\\
 \ds  \frac{-r +\sqrt{ r^2 -\frac {2 \kappa \lambda}{\dz} \left( r\left(\frac \mu \lambda - 
\frac {\phi(k_{\min})}{\kappa \lambda}\right) -\frac {\kappa\lambda } {2 \dz}  (q_{0,j+1})_+ ^2  + g(k_{\min})\right)  }}
{\frac {\kappa\lambda}{\dz} }      
  &\ds \hbox{if }\quad Q>\left(q_{0,j-1}  \right)_- .
    \end{array}
\right.
\end{equation}
Then $p$ satisfies the constraint in  (\ref{eq:65}) if and only if $p\ge q +  \frac {\mu} \lambda -\frac{\phi(k_{\min})}{\kappa \lambda} $, and $B_j$ is
 computed by maximizing a concave and quadratic function on the set  $p\ge q +  \frac {\mu} \lambda -\frac{\phi(k_{\min})}{\kappa \lambda} $.
\end{remark}

\subsection{The discrete scheme at $i=N$}\label{sec:discrete-scheme-at-1}
For brevity, we do not write the 
numerical scheme corresponding to the boundary condition  at  $k=k_{\max}$,  because the equations (two equations for each value of $j$, $0\le j\le M$,) may be obtained  in exactly the same way as in the previous paragraph.

\subsection{Solving the system of nonlinear equations: a long time approximation}\label{sec:solv-syst-nonl}

The system of equations including (\ref{eq:61})-(\ref{eq:62}) for $0<i<N$ and $0\le j\le M$, and the discrete versions of the boundary conditions at $k=k_{\min}$ and $k= k_{\max}$ described above,
  can be written in an abstract form as follows:
\begin{equation}\label{eq:77}
  \cF( \cU, \cP )=0,
\end{equation}
where $\cF$ is a nonlinear map from $\R^{2(N+1)(M+1)}$ to $\R^{2(N+1)(M+1)}$ such that the Jacobian matrix of $\cF( \cU, \cP )$ has negative diagonal entries.

We  aim at solving the discrete system  (\ref{eq:77}) by a long time approximation involving an explicit scheme. The reason for choosing an explicit scheme lies in the complexity of the boundary conditions. Finding an  implicit or semi-implicit scheme consistent  with the nonlinear boundary  conditions seems challenging.

We fix a time step $\dt>0$.

Setting $\cU^\ell = (U^\ell_{i,j} )_{0\le i\le N, 0\le j\le M}$ 
and $\cP^{\ell} = (p^{\ell}_{i,j} )_{0\le i\le N, 0\le j\le M}$,
we compute the sequence $(\cU^\ell,\cP^\ell)$ by the induction:
\begin{equation}
(\cU^{\ell+1},\cP^{\ell+1})=(\cU^\ell,\cP^\ell)-\dt \cF(\cU^\ell,\cP^\ell),
\end{equation}
and expect that the sequence converges as $\ell\to +\infty$. It the latter case, the limit is a solution of (\ref{eq:77}).

\section{Numerical simulations}\label{sec:numer-simul}

\subsection{The parameters}\label{sec:parameters}
The numerical simulations reported below aim at describing some aspects of  the short term dynamics of the oil market. Some of the parameters $(r,\kappa, \lambda,\mu, \epsilon,q_\circ,c)$ come from the calibration of the model \cite{edmond1} to  prices, CAPEX, OPEX, capacities and production  observed in the last three decades. The value of $k_{\max}-k_{\min}$ is qualitative reasonable guess, since there is no direct observation of the arbitraged storage. Indeed, the real storage includes strategic, operational and arbitraged storage. The other parameters and functions are qualitatively reasonable guesses to obtain the proxy of the investments delay effects and the costs of storage.

We believe that these numerical simulations make it possible to explain the sharp drops   in  prices  and the general dynamics of the oil industry that have been observed in 2015 and 2020.


We take
\begin{displaymath}
  \begin{split}
  b(k,z,p)=  a \left(\frac {k_{\max} -k} {k_{\max} -k_{\min}}\right)^2   -a \left(\frac {k -k_{\min}} {k_{\max} -k_{\min}}\right)^2 &+  \kappa( \lambda p -\mu )\\ & \quad \hbox{far enough from $z=z_{\min}$ and $z=z_{\max}$}, \\
  g(k)=0,\\ 
  \end{split}
\end{displaymath}
with 
\begin{displaymath}
\begin{split}
r=0.1,\\
\epsilon=4.\; 10^{-4},\\
  a= 0.01,\\
\kappa=2 \;10^{-3},\quad  \lambda=0.4,\quad \mu=25,\\
q_\circ=0.42,  \quad \alpha=10^4,\\
c=10,\\
\nu_z=10^{-4},\quad \sigma(k)=0, \hbox{ for all } k.
\end{split}
  \end{displaymath}
We set $k_{\min}=0$, $k_{\max}=0.05$,  $z_{\min}=0.35$ and $z_{\max}=0.75$. We refer to Appendix~\ref{sec:appendix} for a simulation in which the cost of storage $g(k)$ is not zero and penalizes situations in which the storage capacities are close to full.\\
 In order to keep the expression of $b$ simple, we have decided not to write explicitly the perturbations of $b$ near  $z=z_{\min}$ and $z=z_{\max}$, which do not impact the solution in the region of interest. In the same vein, recall that $z_{\min}$ and $z_{\max}$ are  technical bounds on $z$ which are only useful for numerical purposes, i.e. in order to work in a bounded domain: another sensible choice of $z_{\min}$ and $z_{\max}$ would not change the solution in the region of interest.
\\
The mesh parameters are $N=M=200$, and the time step is $\Delta t= 0.00001$.

\subsection{The results}\label{sec:resultss}
\begin{center}
  \includegraphics[width=0.8\linewidth]{{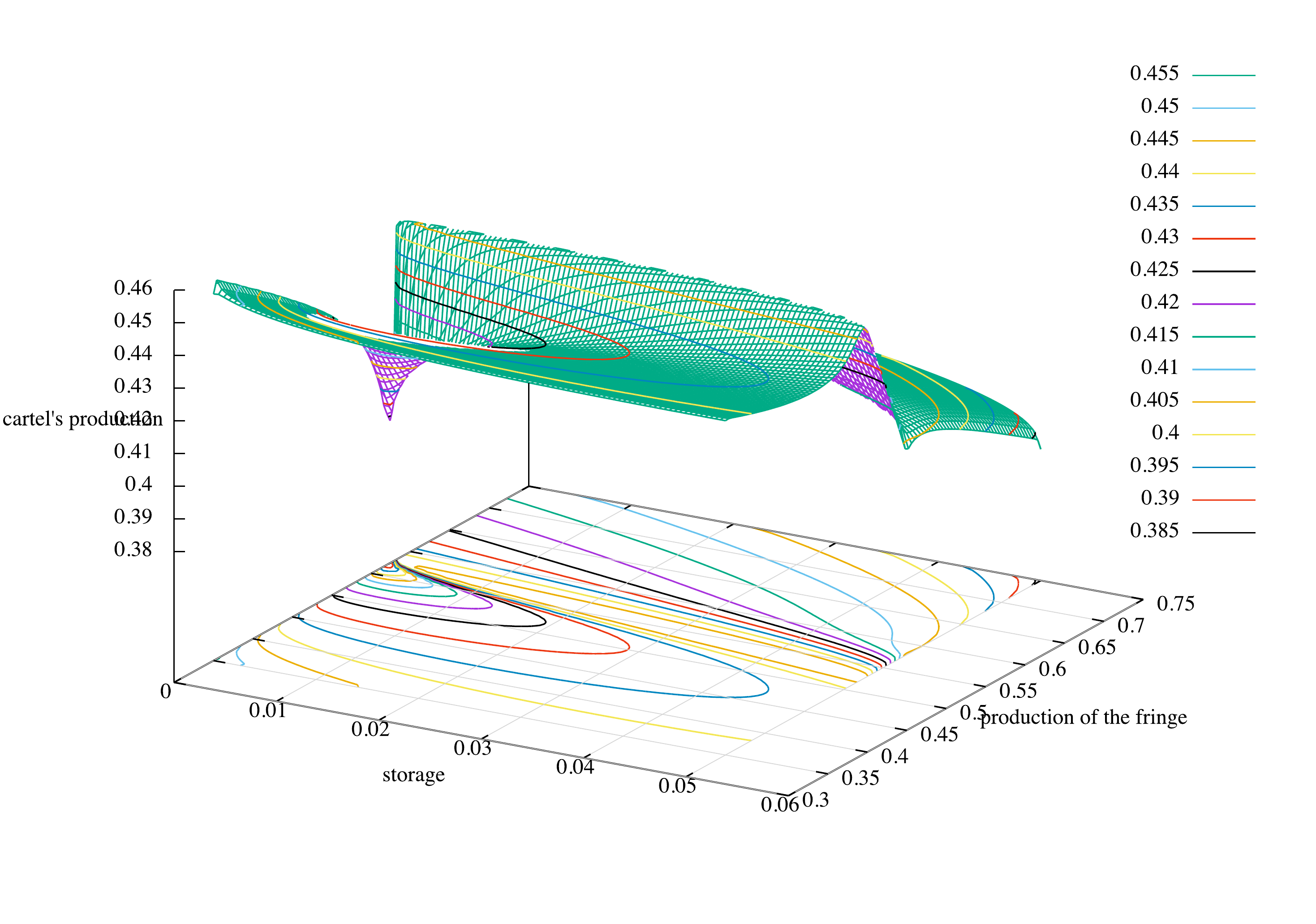}}
  \captionof{figure}{The optimal production level of the cartel as a function of $k$ and $z$: we see that the optimal level of production displays a shock 
whose amplitude is maximal at $k=k_{\min}$ and vanishes at $k=k_{\max}$. }
  \label{fig:1}
\end{center}

\begin{center}
  \includegraphics[width=0.8\linewidth]{{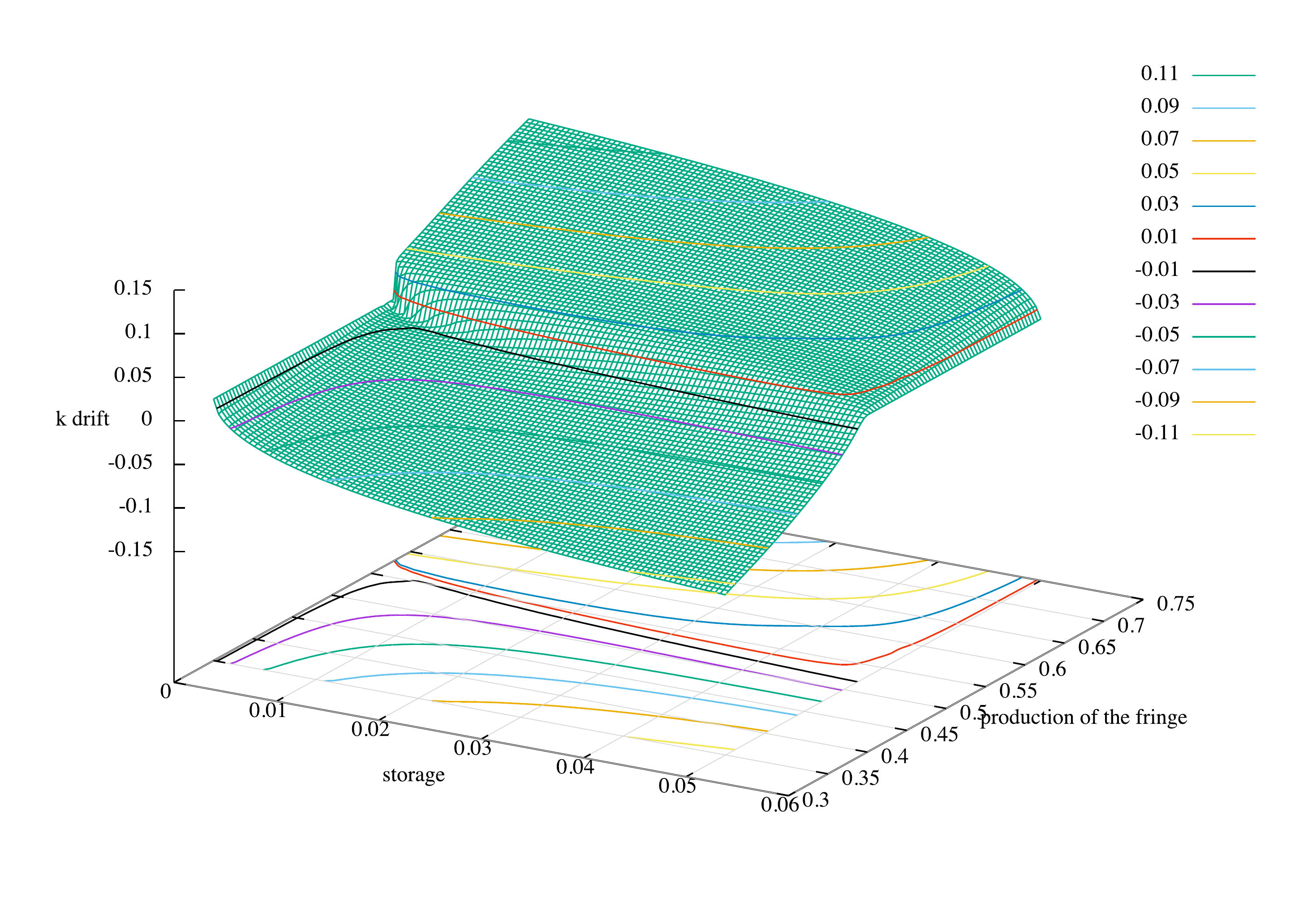}}
  \captionof{figure}{The $k$ component of the optimal drift, which gives the dynamics of the storage level, as a function of $k$ and $z$. Note that it behaves like  $ \sqrt{k-k_{\min}}$  when it is nonpositive near $k=k_{\min}$, and that  it behaves like  $ \sqrt{k_{\max}-k}$  when it is nonnegative near $k=k_{\max}$.}
  \label{fig:2}
\end{center}
  
\begin{center}
  \includegraphics[width=0.8\linewidth]{{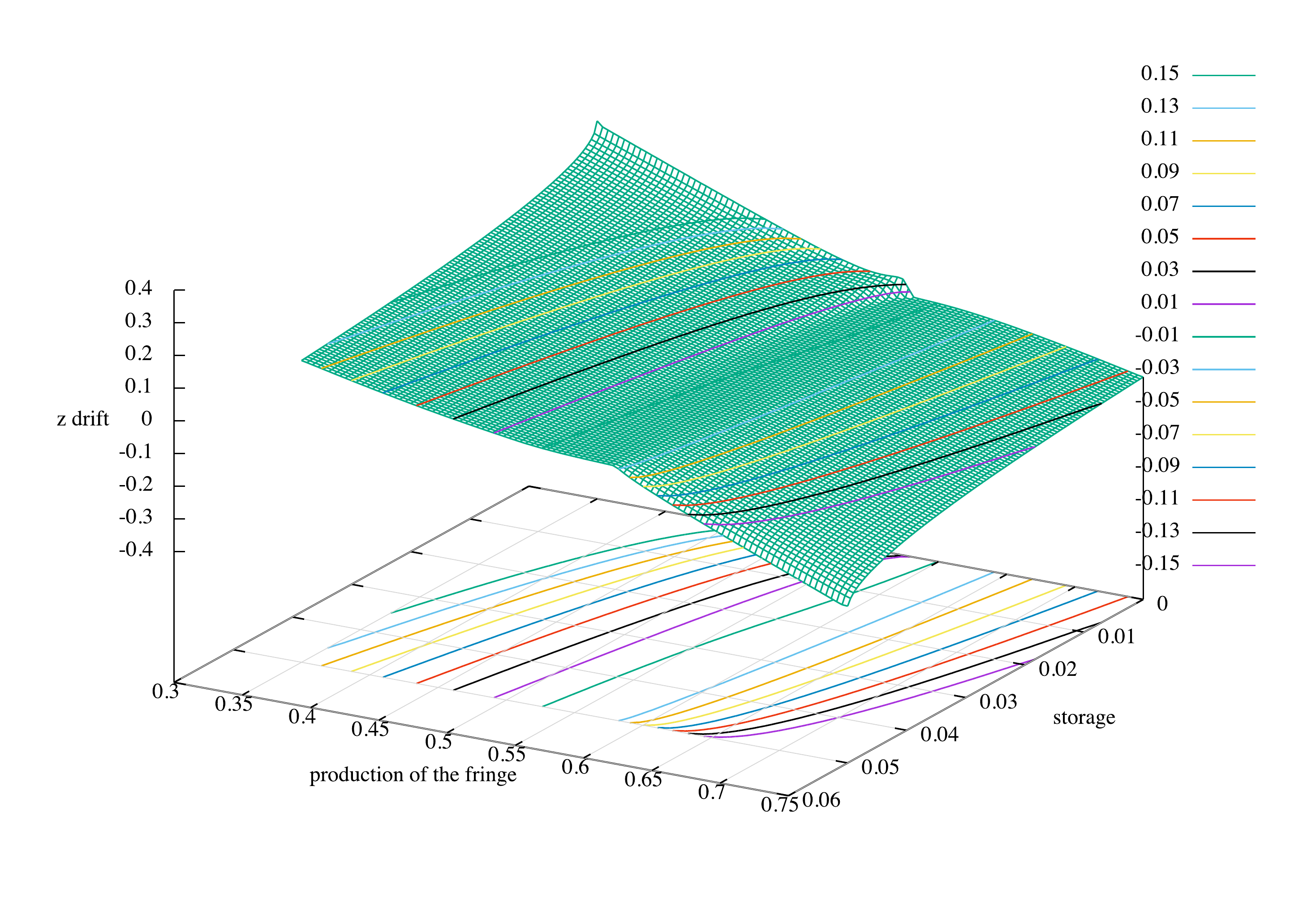}}
  \captionof{figure}{The $z$ component of the optimal drift, which gives the dynamics of the production  level of the fringe, as a function of $k$ and $z$.}
  \label{fig:3}
\end{center}

\begin{center}
  \includegraphics[width=0.8\linewidth]{{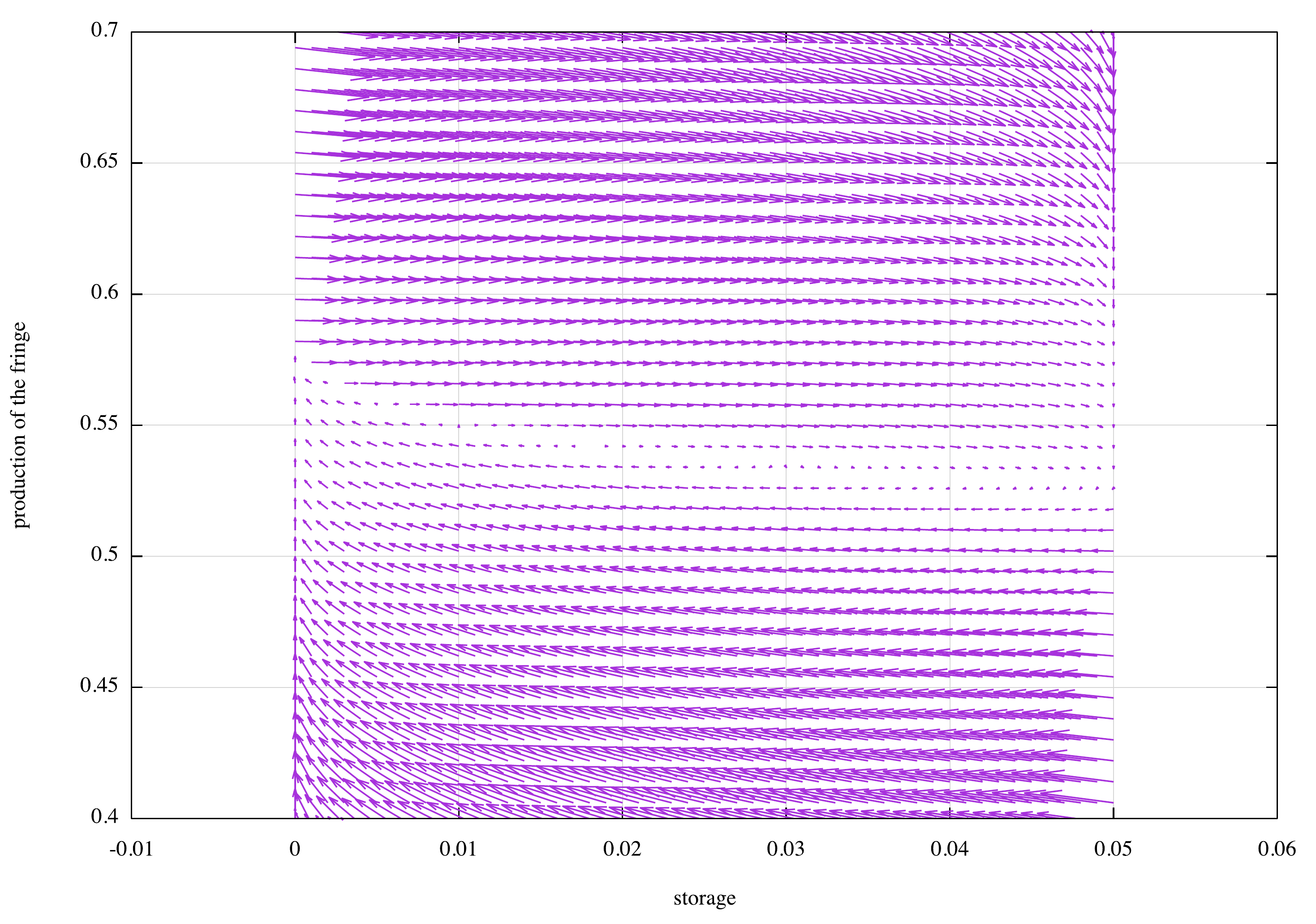}}
  \captionof{figure}{The optimal drift as a function of $k$ and $z$.}
  \label{fig:4}
\end{center}

\begin{center}
  \includegraphics[width=0.8\linewidth]{{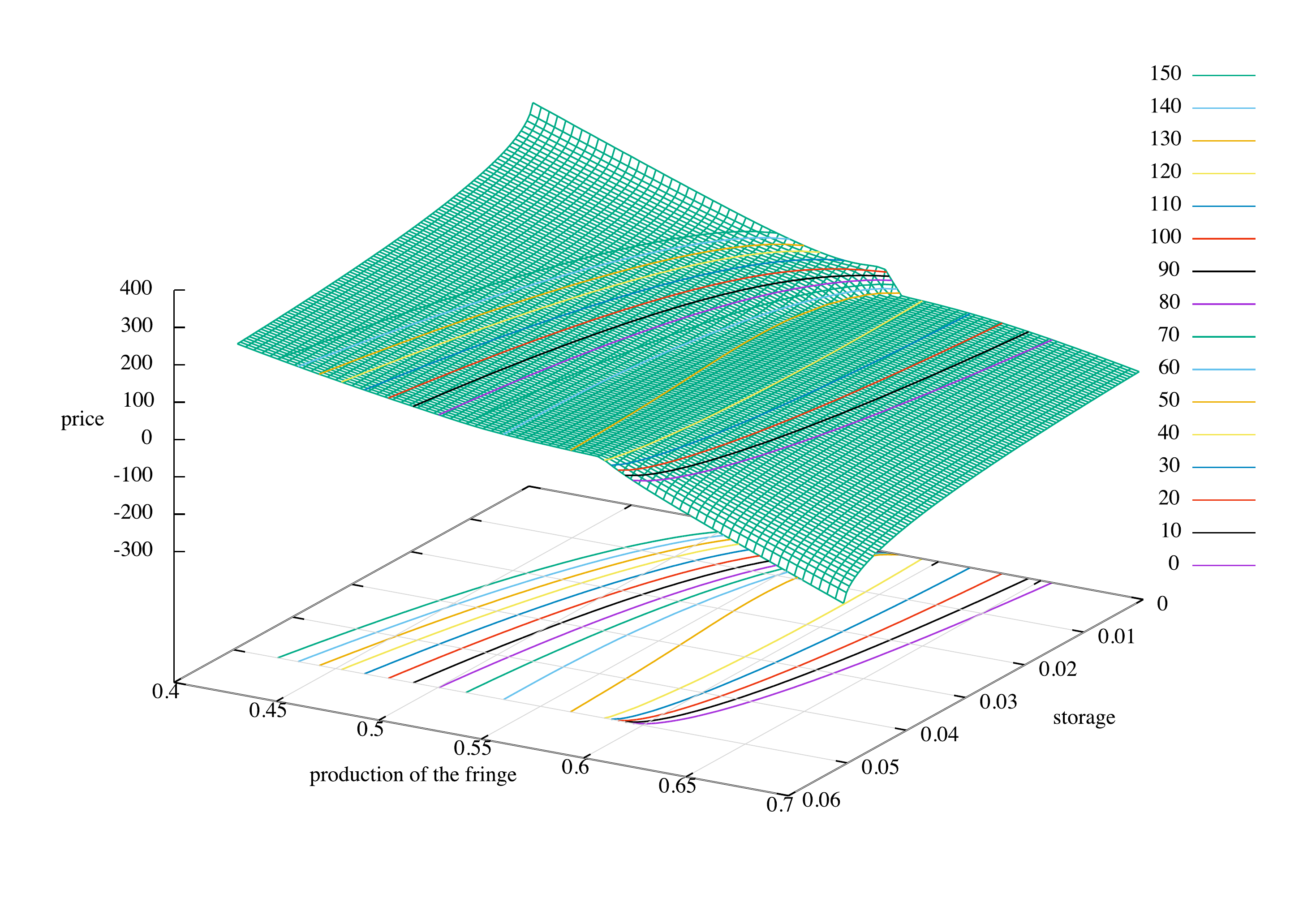}}
  \captionof{figure}{The price as a function of $k$ and $z$;  there is also a shock in the price. Note that the price takes negative values for large values of $z$, but such large values are most often irrelevant in the oil industry (the level of production of the competitive fringe is close to  $0.56$ and does not vary more than  $5\%$). The negative values of $p_t$ are due to the fact that we chose not to put any constraints on the production level $q$.}
  \label{fig:6}
\end{center}

\begin{center}
  \includegraphics[width=0.8\linewidth]{{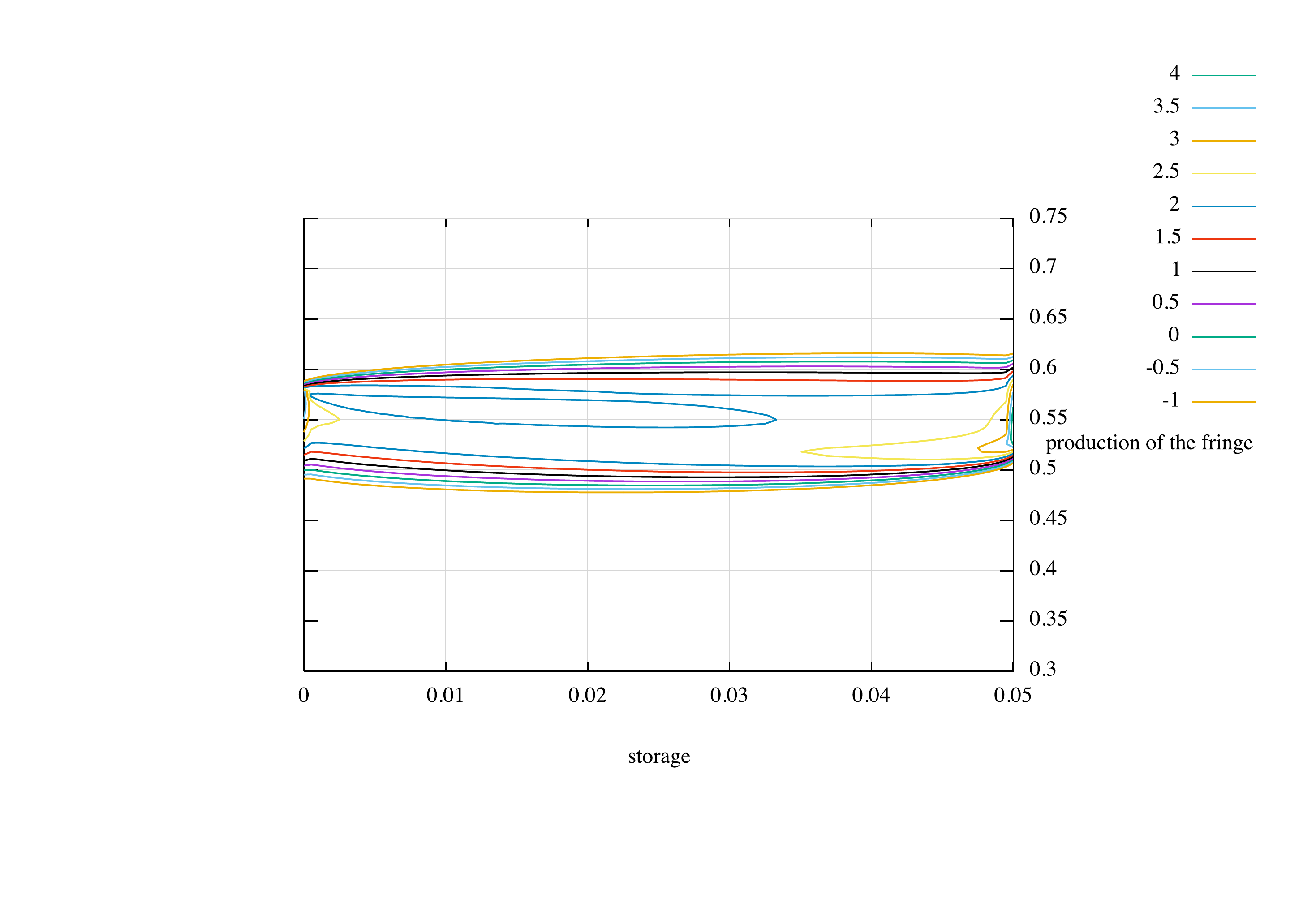}}
  \captionof{figure}{The contours of the invariant measure of $(k_t,z_t)$ (in logarithmic scale):  it is concentrated around the stable  cycle. Within the cycle, the density of the measure is much higher in the  region  close to  the lines $k=k_{\min}$ and $k_{\max}$, because the drift is small there.}
  \label{fig:5}
\end{center}

\begin{center}
  \includegraphics[width=0.45\linewidth]{{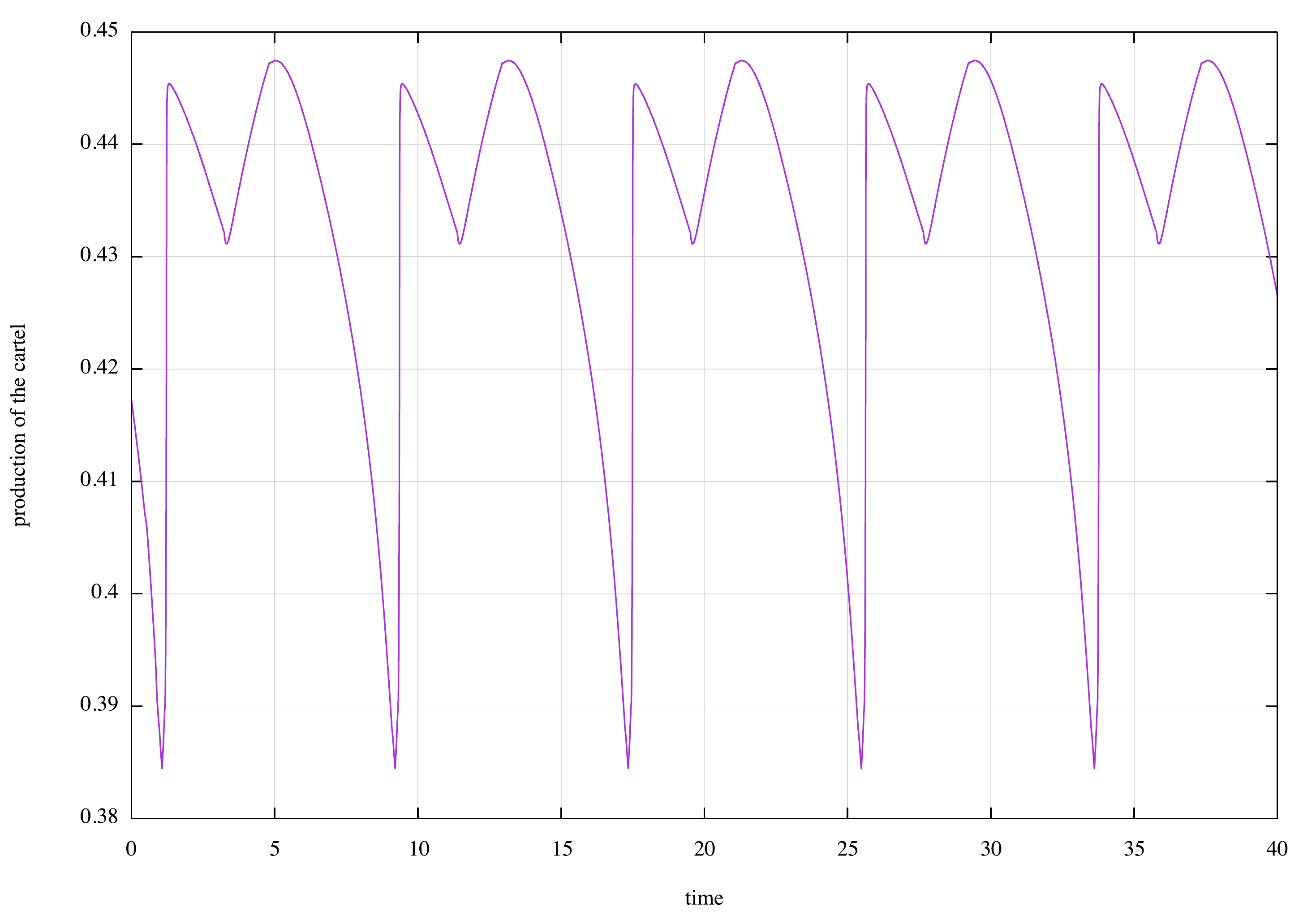}}
\includegraphics[width=0.45\linewidth]{{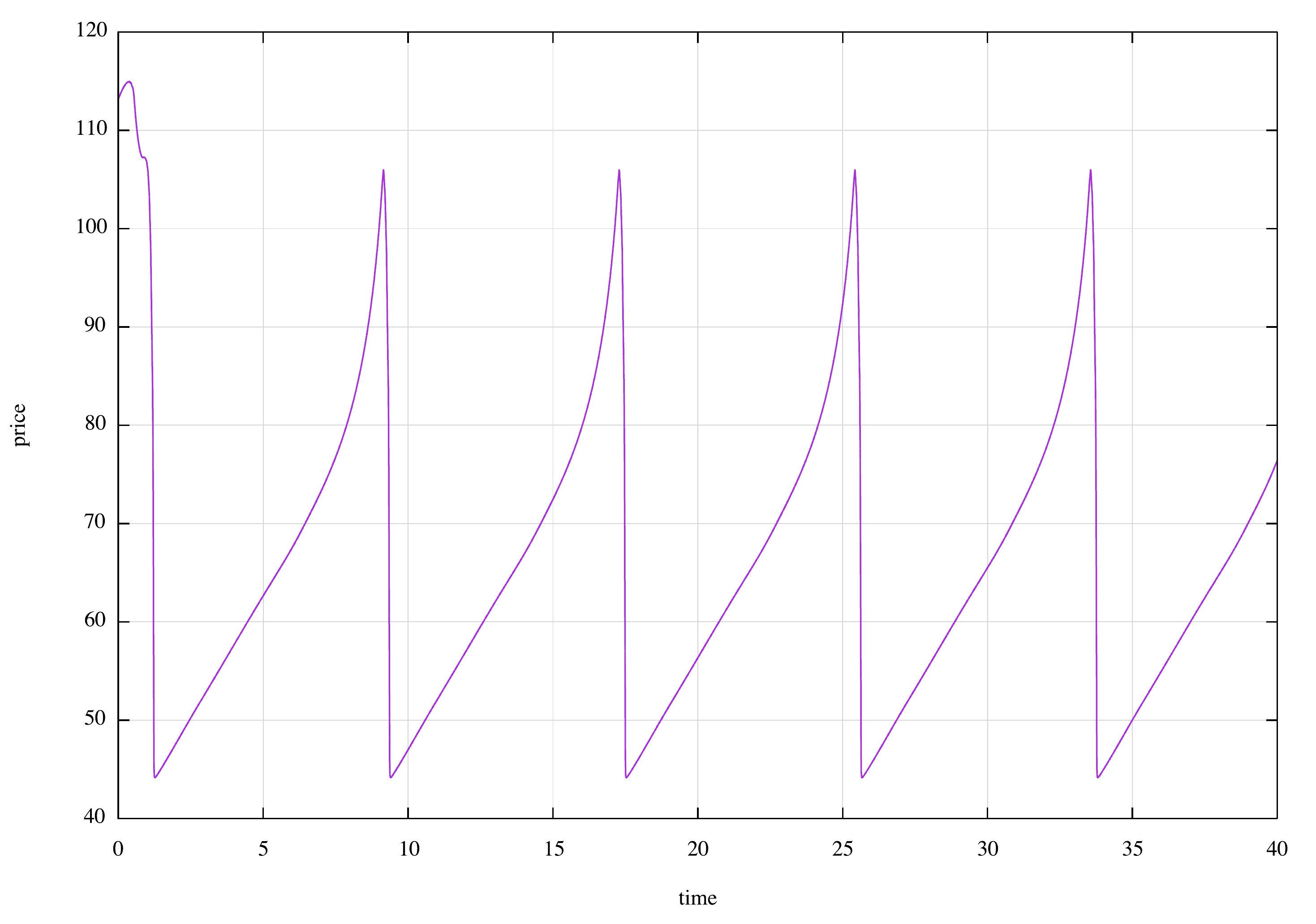}}
\\
\includegraphics[width=0.45\linewidth]{{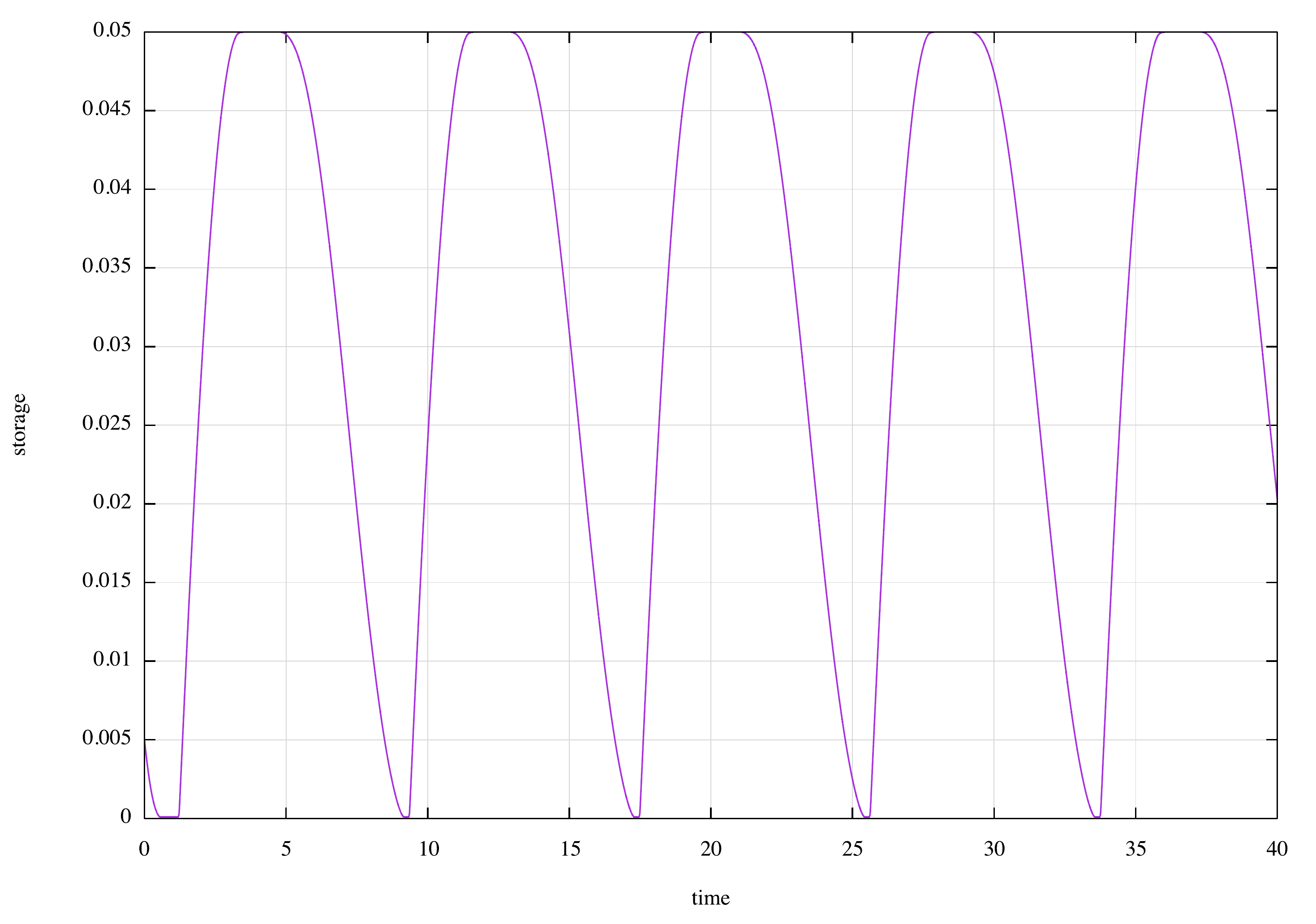}}
\includegraphics[width=0.45\linewidth]{{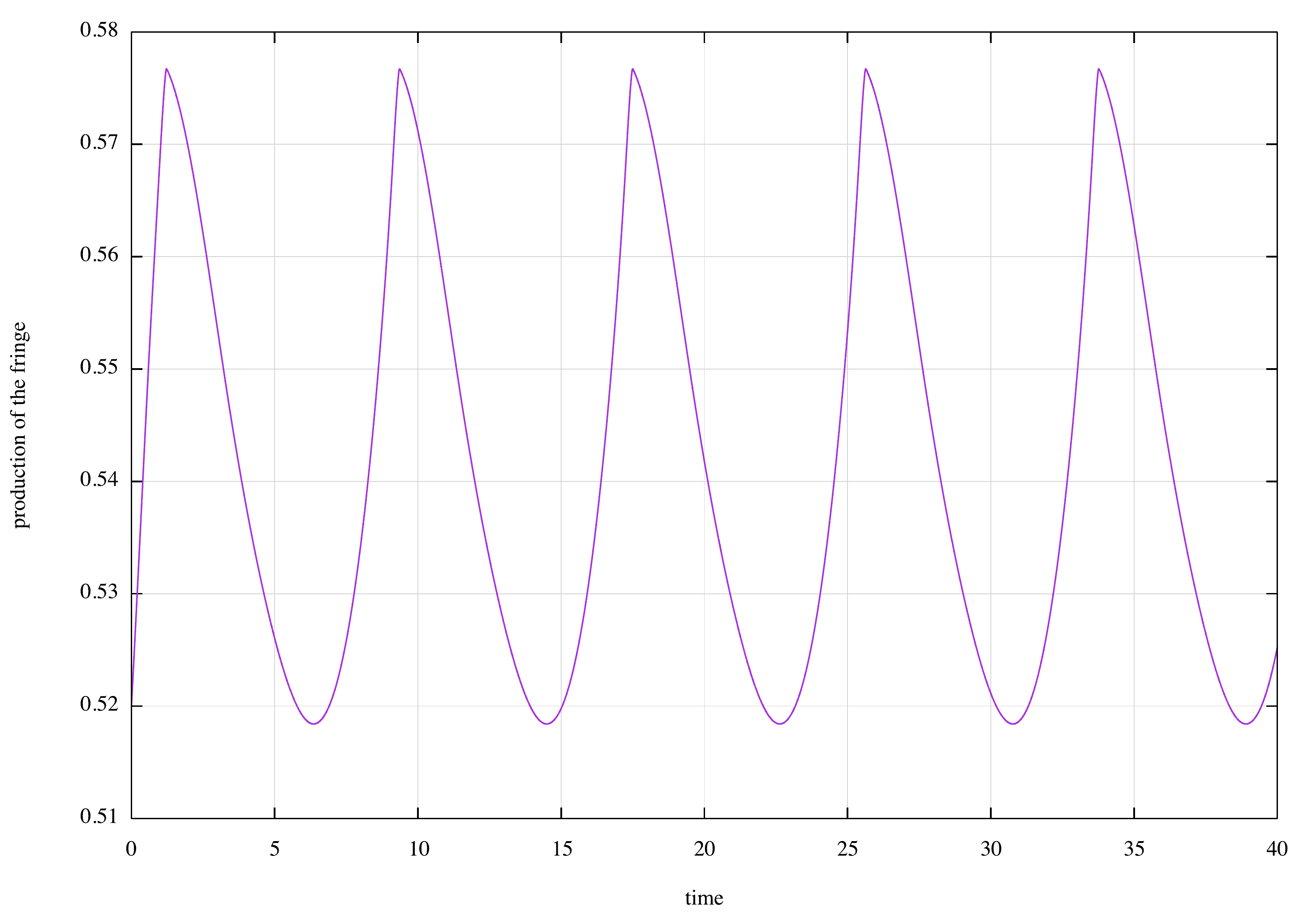}}
  \captionof{figure}{ A stable cycle:   with the functions  $u$ and $p$ computed numerically by the finite difference method,  we neglected the noise (which is besides small in the present case) and simulated a trajectory starting from the point $(k,z)=(0,0.5)$ by means of a standard Euler scheme applied to 
(\ref{eq:dynz})-(\ref{eq:3}) (the Gaussian noise on $z$ is therefore not taken into account): we see that after a small time, the trajectory becomes periodic with respect to time, with a period of the order of $7.5$ years. In reality, there is noise, hence the cycles are not so evident and so regular. Top-Left: the production of the cartel vs. time. Top-Right: the price vs. time. Bottom-Left: the  level of storage vs. time. Bottom-Right: the production of the fringe vs time.}  \label{fig:cycles}
\end{center}

\paragraph{Comments }$\;$

The numerical results are displayed on Figures  \ref{fig:1}-\ref{fig:cycles}:
Figure  \ref{fig:1} contains the graph and the contour lines of the optimal  production level $q^*$ of the cartel as a function of $k$ and $z$:   $q^*$ is discontinuous across a line in the $(k,z)$ plane (the discontinuity is smeared a little due to the small diffusion). The amplitude of the discontinuity is maximal  at $k=k_{\min}$ and vanishes at $k=k_{\max}$, i.e. at $k=k_{\max}$, the optimal production $q^*$ depends continuously on $z$.  Further comments on the optimal policy will be made below. Figures  \ref{fig:2}-\ref{fig:4} give information on the dynamics of $k$ and $z$ resulting from the optimal policy. Figure  \ref{fig:2} contains the graph and the contour lines of the drift of the storage level $k$, i.e. $q^* +z -D(p^*)$, as a function of $k$ and $z$. It is of course discontinuous  across  the same shock line as $q^*$. Roughly speaking, the region located above (respectively below) the shock in the $(k, z)$ plane corresponds to a regime when  storage increases, (respectively decreases).  Note the characteristic behaviour of the optimal drift of $k$  for small values of $z$ near $k=k_{\min}$: it is negative and  vanishes at $k=k_{\min}$ like $\sqrt{k-k_{\min}}$, in agreement with the theoretical results contained in paragraph \ref{sec:few-results-math}. The same behavior is observed for large values of $z$ near $k=k_{\max}$: the optimal drift is positive and  vanishes at $k=k_{\max}$ like $\sqrt{k_{\max}-k}$.
The $z$-component of the drift, i.e. $ b(k, z, p^* )$, discontinuous  across the same shock line, is displayed on  Figure  \ref{fig:3}.
the optimal drift vector in the $(k,z)$ plane is plotted on  Figure  \ref{fig:4}.
 Figure  \ref{fig:6}  contains the graph and the contour lines of $p^*$ as a function of $k$ and $z$, which also has  a shock.
On Figure  \ref{fig:5}, we display the contour lines of the invariant measure of the process $(k_t,z_t)$: we see that it is concentrated around a cycle that takes place around the shock line.
Finally, in Figure \ref{fig:cycles},  with the functions  $u$ and $p$ computed numerically by the finite difference method,  we neglected the noise and simulated a trajectory starting from the point $(k,z)=(0,0.5)$ by means of a standard Euler scheme: we see that after a small delay, the trajectory becomes periodic with respect to time, with a period of the order of $7.5$ years. 
 This will be commented and interpreted below. Note that we purposely initialized the trajectory at $(k,z)=(0,0.5)$ which does not belong to the cycle, in order to illustrate the fact  that the cycle is attractive.

\bigskip

It is remarkable that the system of PDEs   \eqref{eq:6}-\eqref{eq:7}, supplemented with the boundary conditions linked with the constraints on the state variable $k$ and discussed in paragraph \ref{sec:bound-cond-at-1}, leads to a discontinuous optimal strategy   as in Figure \ref{fig:1}. The discontinuous solutions obtained in the simulations may seem surprising. We are going to explain why, on the contrary, the simulated  optimal policy matches qualitatively  well what has been observed in the past few years. We are first going to see that a cyclic strategy of the cartel is naturally linked to the singularity displayed on  Figure \ref{fig:1}. After having described the cycle, we will explain how these results shed a light on what has been observed in 2015 and in 2020.

Note that the level of noise in the present simulation is smaller than it is actually. We have underestimated the  noise in order to shed some light on the mechanism resulting from the model, which would appear less clearly in the presence of noise. Note also that,  knowing $u$ and $p$ from the simulation of the system of PDEs~(\ref{eq:6})-(\ref{eq:7}),
we purposely simulated (\ref{eq:dynz}) and (\ref{eq:3}) instead of (\ref{eq:92}) and (\ref{eq:3}), (recall that $\sigma=0$), in order
  to exhibit a stable cycle  and its  time periodicity. 
As we can see on Figures  \ref{fig:4}, \ref{fig:5} and \ref{fig:cycles}, the present model, in the absence of noise,  leads to a cyclic behaviour  on the time scale of few years (more details will be given below).  Since there are important noise and risk factors in the oil market, such a cyclic behaviour is not always observed, but is has actually been observed twice in the recent past, in 2015 and in 2020.

\subsection{Interpretation of the observed cycle}\label{sec:interpr-observ-cycle}
The cycle that is observed in the numerical simulations, see Figures \ref{fig:5} and \ref{fig:cycles}, is drawn schematically on Figure \ref{fig:0}.  In the absence of substantial randomness, we can  make out four phases.
\begin{center}
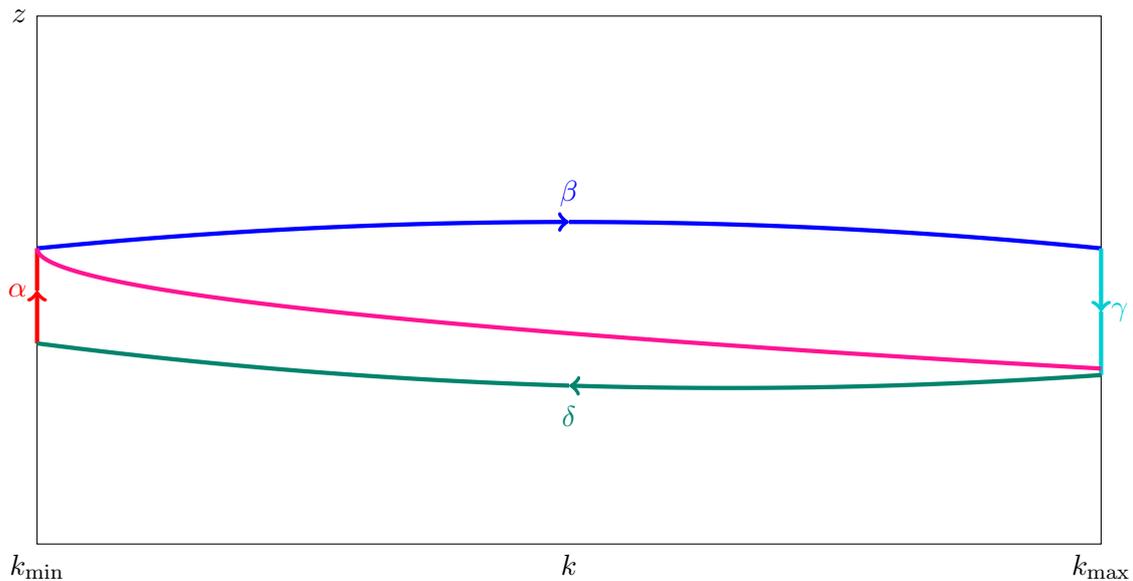

  \begin{tikzpicture}[scale=1.4]
    \draw  (0,0) rectangle (10,5);
    \draw(0,0)   node[below] {$k_{\min}$ };
    \draw(10,0)   node[below] {$k_{\max}$ };
    \draw(5,0)   node[below] {$k$ };
    \draw(0,5)   node[left] {$z$ };
    \draw[->, color=red, ultra thick]   (0,1.9) -- (0,2.4)  ;
     \draw[ color=red, ultra thick]   (0,2.4) -- (0,2.8)  ;
    \draw(0,2.4) node[left, color=red] {$\alpha$};

    \draw[->, color=blue, ultra thick, domain=0:5] plot[id=0112+] function{ 2.8+  0.01*x*(10-x) } ;
    \draw[ color=blue, ultra thick, domain=5:10] plot[id=0112++] function{ 2.8+  0.01*x*(10-x) } ;
    \draw(5,3.1) node[above, color=blue] {$\beta$};   
    \draw[->, color=mydarkturquoise, ultra thick]   (10,2.8) -- (10,2.2)  ;
     \draw[ color=mydarkturquoise, ultra thick]   (10,2.2) -- (10,1.6)  ;
     \draw(10,2.2) node[right,  color=mydarkturquoise] {$\gamma$};
     
     \draw[ <-, color= PineGreen, ultra thick, domain=5:10] plot[id=0113+] function{ 1.6 - 0.01*(x-3)*(10-x) } ;
      \draw[color= PineGreen, ultra thick, domain=0:5] plot[id=0113++] function{ 1.6 - 0.01*(x-3)*(10-x) } ;
    \draw(5,1.4) node[below, color= PineGreen] {$\delta$};

    \draw[ color= mydeeppink, ultra thick, domain=0:10] plot[id=0114+, samples=1000] function{ 2.8 - 0.36*sqrt(x) } ;
    
   \end{tikzpicture}
  \captionof{figure}{\label{fig:0} The cycle $\alpha, \beta, \gamma, \delta$ in the $(k, z)$ plane. In pink, the shock line.}
\end{center}

Below, we describe the four phases and relate them to Figure~\ref{fig:cycles}, focusing for example on the time period  approximately from  $9.5$ to $17$ (years).
\begin{description}
\item[Phase ($\alpha$):] in this phase, which corresponds to the time interval  approximately from  $16.5$ to $17$  (or $9$ to $9.5$) on Figure~\ref{fig:cycles},  storage is close to  minimal, i.e. $k=k_{\min}$
  (recall that we only deal with the storage managed by the arbitrageurs). The monopolistic cartel has therefore the power to drive the price up by  maintaining a low level of production.   When the price goes up, the fringe producers invest in new production capacities and increase gradually their market share. At some point,  the monopolistic cartel would like to  drive  price down. However, the monopolistic cartel is aware that when it  does so by increasing its own production,  arbitrageurs start storing the resource, which  diminishes the intended impact of the policy.   In order to prevent the period of low prices from lasting too long, it is therefore  optimal for  the monopolistic cartel to brutally increase its production, thereby initiating the phases $\beta, \gamma, \delta$ that  bring the monopolistic cartel back to its zone of profit. This is what shows the numerical simulation of Edmond 2; the strategic shock appears clearly on Figure  \ref{fig:1}: indeed, it can be seen that when  storage is minimal, the production of the monopolistic cartel increases brutally when $z$, the production of the fringe crosses a critical value. Then, on Figure \ref{fig:6}, we see that this brutal increase  in production makes the price fall. This phenomenon has been observed in 2015 and 2020.
\item[Phase ($\beta$):] the time interval  approximately from  $9.5$ to $11.5$ on Figure~\ref{fig:cycles}. When the price has fallen, the stored quantity of resource  increases rapidly, and the state $(k_t, z_t)$  is drifted to the right with a velocity nearly parallel to the $k$-axis.  After having  brutally increased its  production, the cartel may  let it decrease smoothly until  storage gets full. In this regime, the arbitrageurs fix the price and  the price increases almost linearly with respect to time, at a rate close to the interest rate $r=10\%$.
\item[Phase ($\gamma$):] the time interval  approximately from  $11.5$ to $13$ on Figure~\ref{fig:cycles}. Storage is full. The monopolistic cartel can now increase its production again and maintain the low level of price as long as necessary in order to deter the fringe from investing   or even to diminish its existing capacities (see below). The production of the fringe, $z$, decreases to the value that suits the cartel.
\item[Phase ($\delta$):] the time interval  approximately from  $13$ to $16.5$ on Figure~\ref{fig:cycles}. Since the value of $z$
is low enough, the  monopolistic cartel may reduce its production. Then the stored quantity of resource decreases and the price is driven up more rapidly than in the phases $ \beta$ and $\gamma$.
   When  storage is empty, the cartel can start the phase $\alpha$ of the cycle and raise  price to the optimal value.
\end{description}
Even if some aspects may vary with the choice of parameters, see for example Appendix~\ref{sec:appendix}, the main features discussed above seem robust with respect to the choice of the parameters: the discontinuity and the cycle comprising a low price period in order for the cartel to recover market shares and a high price period leading to large profits. One can also see that the present model gives rise  to backwardization and contango periods: indeed, taking $p_{t+1}-p_t$ as an approximation for the slope of the futures curve, Figure ~\ref{fig:cycles} shows that backwardization (resp. contango) occurs when  storage is empty (resp. full). Note that given the noise, backwardization and contango may also occur when the constraints on  storage are still not binding.

\subsection{Discussion on what happened in 2015 and 2020}\label{sec:disc-what-happ}
\begin{description}
\item[ In 2015,]  OPEC had decided to reconquer the  market share  that had been lost due to the fast development of the US shale industry from 2009 to 2015. The price drop was then strong and sudden: prices dropped from  $ \$ 100$  to $ \$ 40 $ per barrel in  a few months. At that time, this price  drop  was analyzed as an attack against US shale. In the spirit of the present model,  it was rather an attack against all  competitors, aimed at recovering  market share. Indeed, the OPEC strategy had a strong impact not only on the US shale industry,  which in fact, proved  strong  resiliency and coping abilities, but in many other fringe producers.
\item[ In 2020,] things happened in a less ``classical'' way. A ``rare disaster'' occured.  An exogenous shock to demand occurred in the first quarter of the year, with a magnitude of the order of ten times the standard deviation of the usual shocks to demand . Although it was not designed to handle such situations, our model seems to give a  good explanation of what happened.
In order to understand 2020, let us recall that in the present model, $z$ (resp. $q$) is the ratio of the non OPEC  capacity of production to the global level of demand, (resp. the ratio of the OPEC production to the global level of demand). Before 2020, the global level of demand increased regularly from year to year (and rather slowly) with an annual  growth rate of the order of $2\% \pm 1 \%$. In the first semester of 2020, the sanitary crisis resulted in a sudden, unexpected and exceptional drop of the global demand, of the order of $10\%$ to $15\%$. Therefore, the variable $z$ got suddenly increased by  $10\%$ to $15\%$, and the monopolistic cartel got 
carried to the upper side of the shock. In our model, the optimal response was to increase immediately production. This is precisely what  happened to the surprise of many. Many observers have considered that what seemed a conflict between OPEC and Russia, which led to an increase in production while the demand collapsed, was suicidal.  However, from the viewpoint of our model, this strategy was simply intended at entering the phases $\beta,\gamma, \delta$ of the cycle, which drive $z$ to the value desired by the cartel. In other words, the aim was to reduce the capacity of production of the competitors as fast as possible. Indeed, as soon as OPEC had increased its production after the collapse of the demand,  storage became rapidly full (stage $\beta$). The maximal level $k_{\max}$ was reached in a few weeks, and the cartel could drive prices to a very low level (prices even went  negative  during a very short period). Many planned investments  stopped, and some production units were  definitely closed. After this period, OPEC started to strongly reduce its production. The production drop was strong in absolute value, but not so strong relatively to the global  demand, hence in  good agreement with our model (recall that all the quantities in the present model are reduced by taking the ratio over the global demand). Then, in the second semester of 2020 and in 2021, the demand should increase, which would imply a fast decrease of $z$.  \\
Hence, despite the fact that the collapse of the global demand in 2020 was  rare event, (seven to ten times the standard deviation of the historical shocks to global demand), our model gives a satisfactory qualitative explanation of the cycle $\alpha, \beta, \gamma, \delta$ that is being observed (in a very accelerated version) in 2020.
\end{description}

\section{Conclusion}
\label{sec:conclusion}
We have proposed a model  which presents new aspects:
\begin{itemize}
\item new ideas on the interactions between a cartel (dominant player) and a crowd of arbitrageurs, 
emphasizing the impact of the constraints on the arbitrageurs,
\item a new system of coupled non linear PDEs and boundary conditions,
\item a new kind of  discontinuous optimal strategy  for the dominant player facing the crowd, 
\item original numerical schemes, robust enough to catch the latter discontinuous solution.
\end{itemize}
Finally, the output of the model shed  original light on an unprecedented crisis.
\begin{acknowledgement}
  
  This research was supported by Kayrros.
   Y. Achdou, C. Bertucci, J-M. Lasry and P-L. Lions  acknowledge partial support from the Chair Finance and Sustainable Development and the FiME Lab  (Institut Europlace de Finance).
  Y Achdou and C. Bertucci    acknowledge partial support from the ANR
  (Agence Nationale de la Recherche) through MFG project ANR-16-CE40-0015-01.
\end{acknowledgement}
\bibliographystyle{plain}
\bibliography{edmond_2}

\appendix \section{A simulation with a cost of storage}\label{sec:appendix}
Here, we keep the parameters as in Section~\ref{sec:numer-simul}, except that we take $k_{\max}=0.07$ and we suppose that there is a cost of storage, which has the form 
\begin{displaymath}
  g(k)=10 \left(\frac { k- k_{\min}    } { k_{\max}- k_{\min}    } \right)^3.   
\end{displaymath}
Such a cost heavily penalizes the situations in which  storage is close to full.  \\
The results have been obtained using exactly the same method as in Section~\ref{sec:numer-simul}, and are represented on Figures  \ref{fig:21},  \ref{fig:22},  \ref{fig:24},  \ref{fig:26},  \ref{fig:25} and  \ref{fig:2cycles}. The results look rather similar, but since $k_{\max}$ is larger than in  Section~\ref{sec:numer-simul}, and since that near to full storage is penalized, these is no phase $\gamma$ in the cycle, as can be seen on Figure  \ref{fig:2cycles}.

\begin{center}
  \includegraphics[width=0.8\linewidth]{{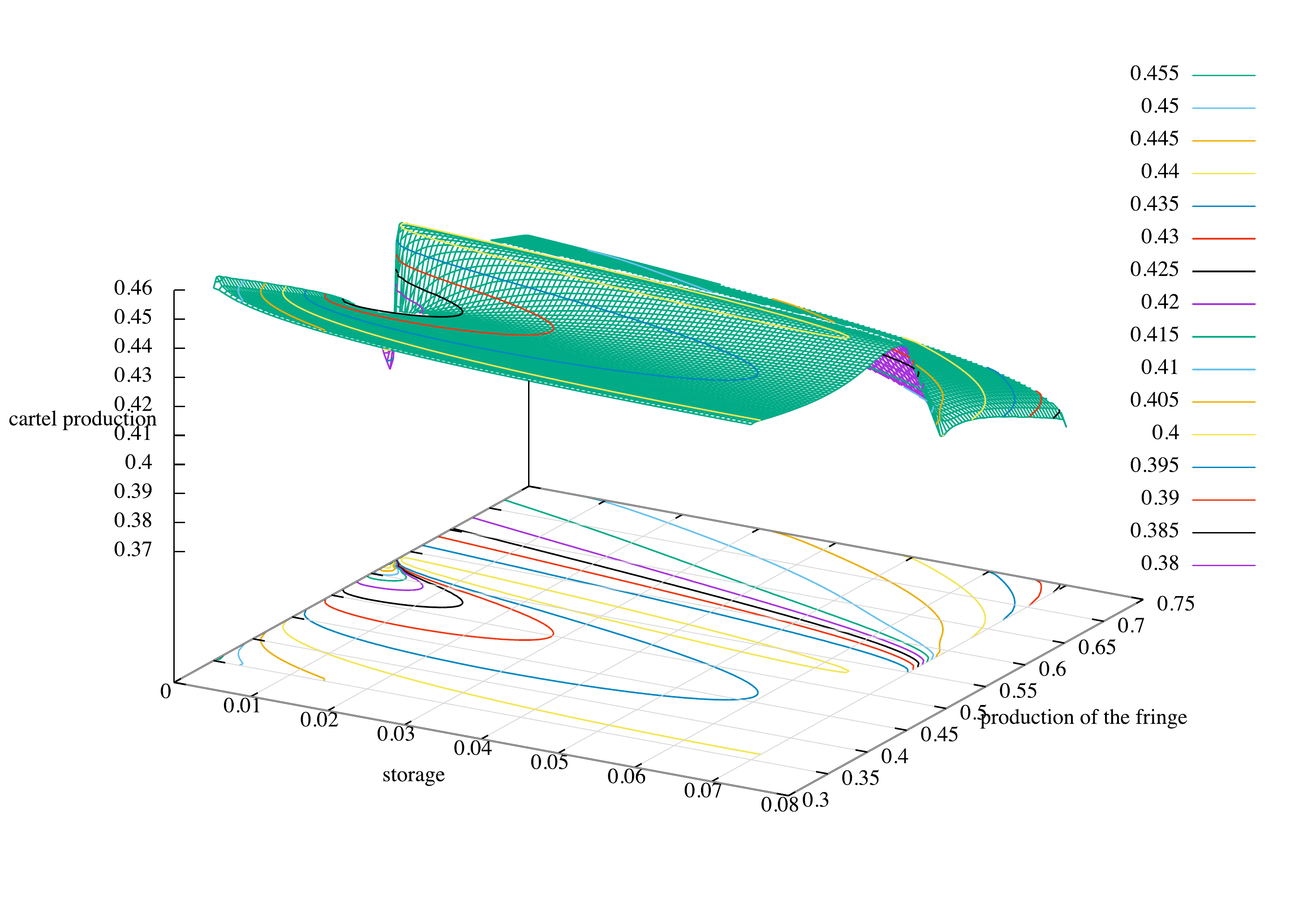}}
  \captionof{figure}{Case 2: the optimal production level of the cartel as a function of $k$ and $z$}
  \label{fig:21}
\end{center}

\begin{center}
  \includegraphics[width=0.8\linewidth]{{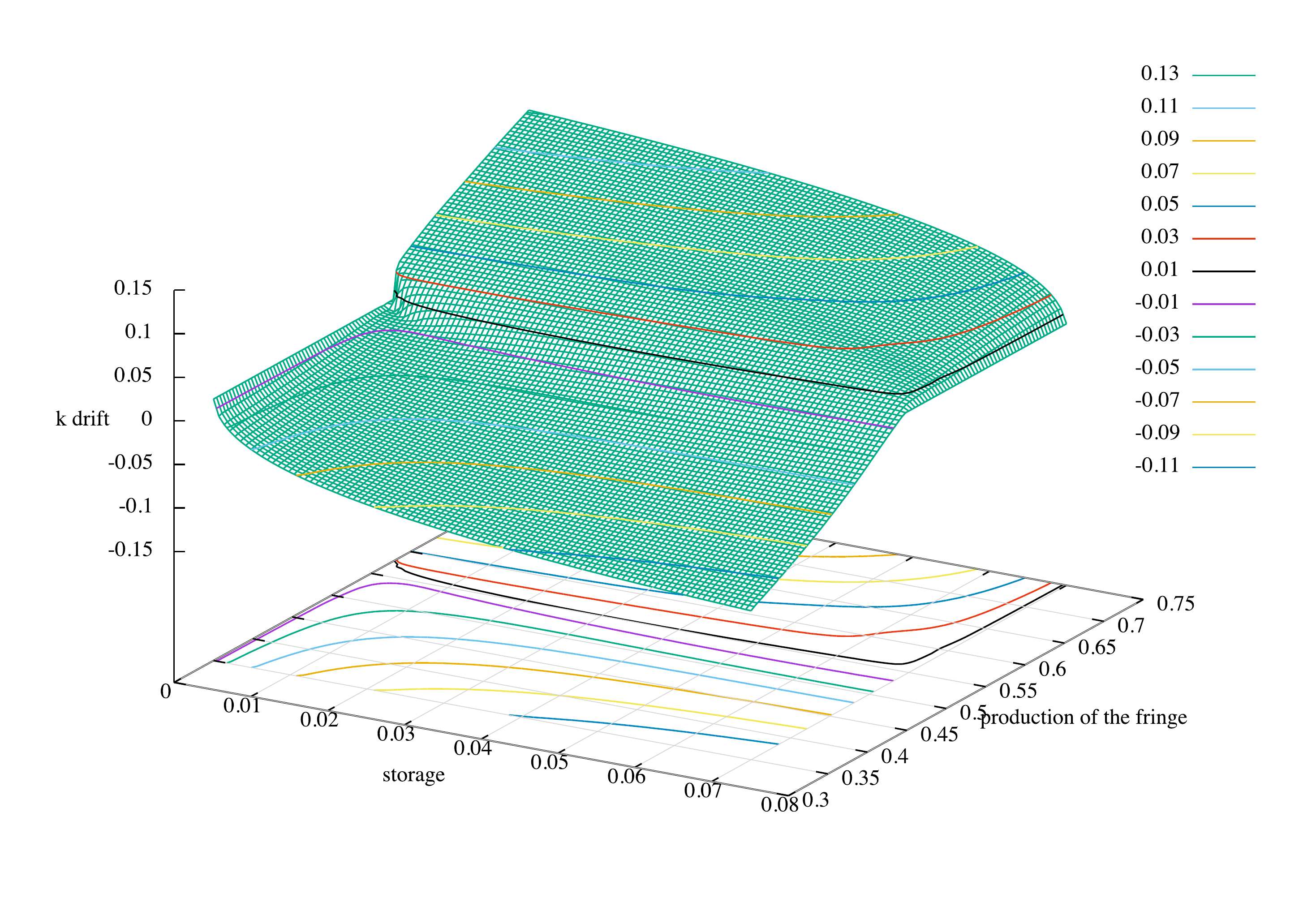}}
  \captionof{figure}{Case 2: the $k$ component of the optimal drift, which gives the dynamics of the storage level, as a function of $k$ and $z$. }
  \label{fig:22}
\end{center}

\begin{center}
  \includegraphics[width=0.8\linewidth]{{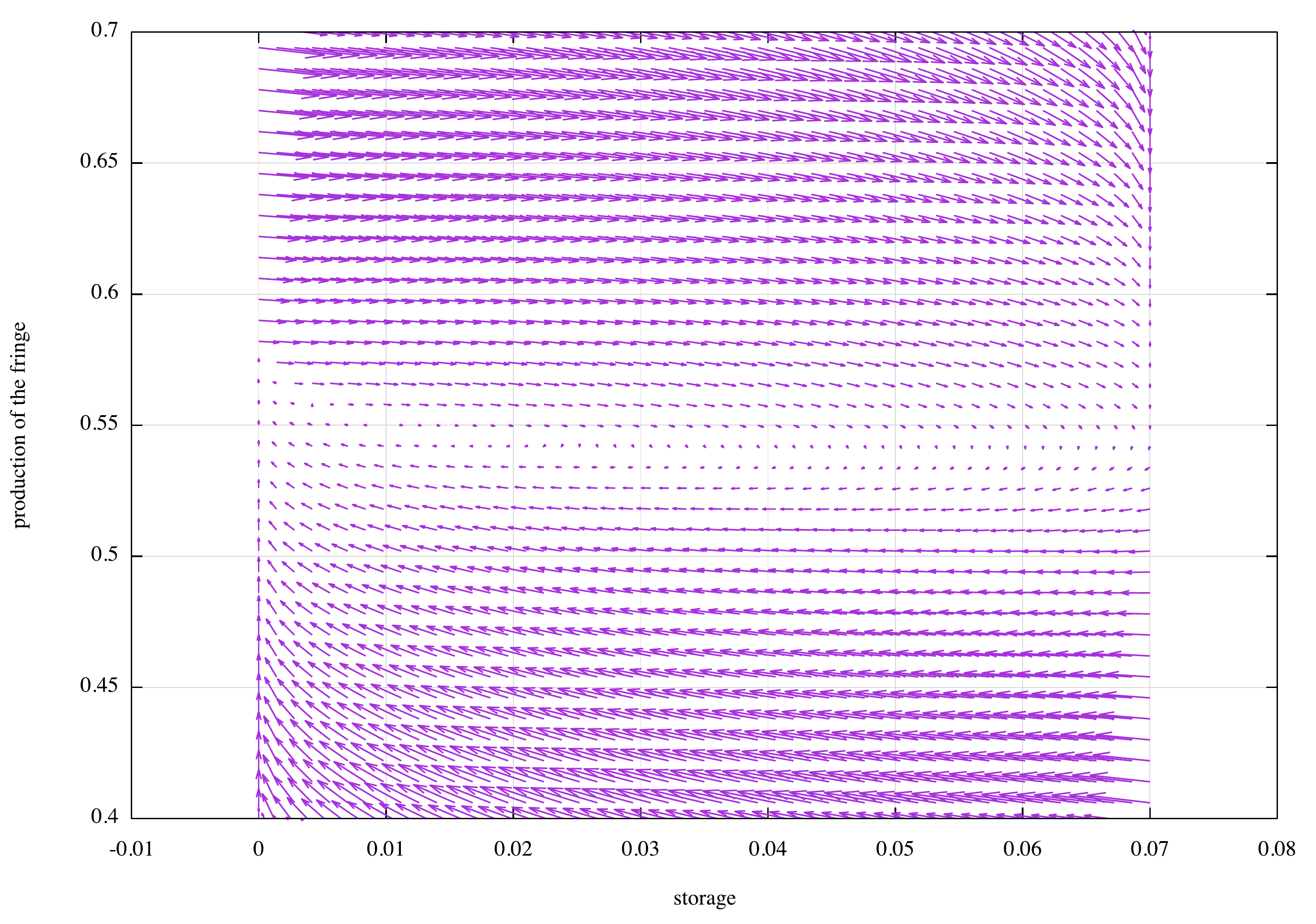}}
  \captionof{figure}{Case 2:  the optimal drift as a function of $k$ and $z$.}
  \label{fig:24}
\end{center}

\begin{center}
  \includegraphics[width=0.8\linewidth]{{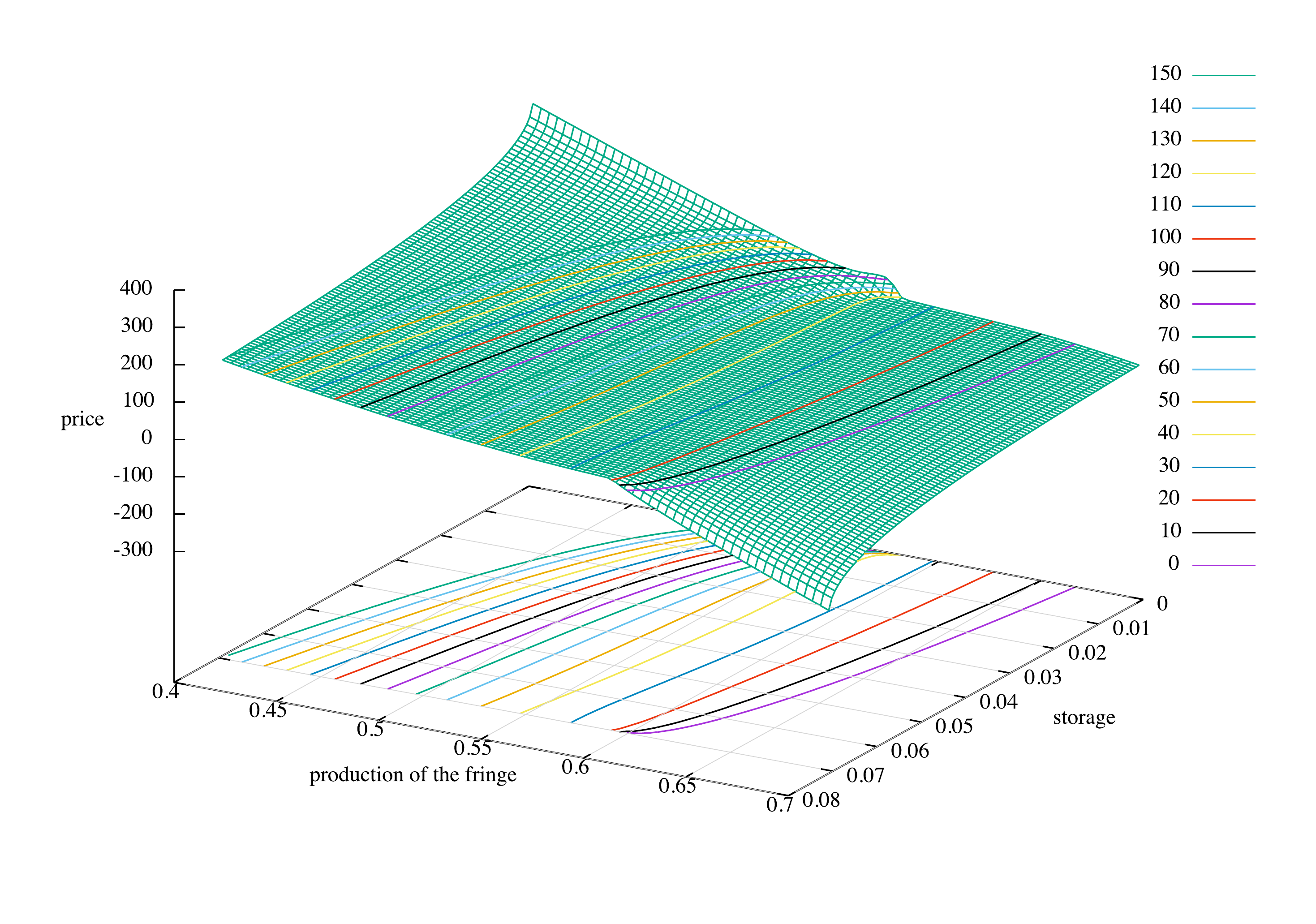}}
  \captionof{figure}{Case 2: rhe price as a function of $k$ and $z$.}
  \label{fig:26}
\end{center}

\begin{center}
  \includegraphics[width=0.8\linewidth]{{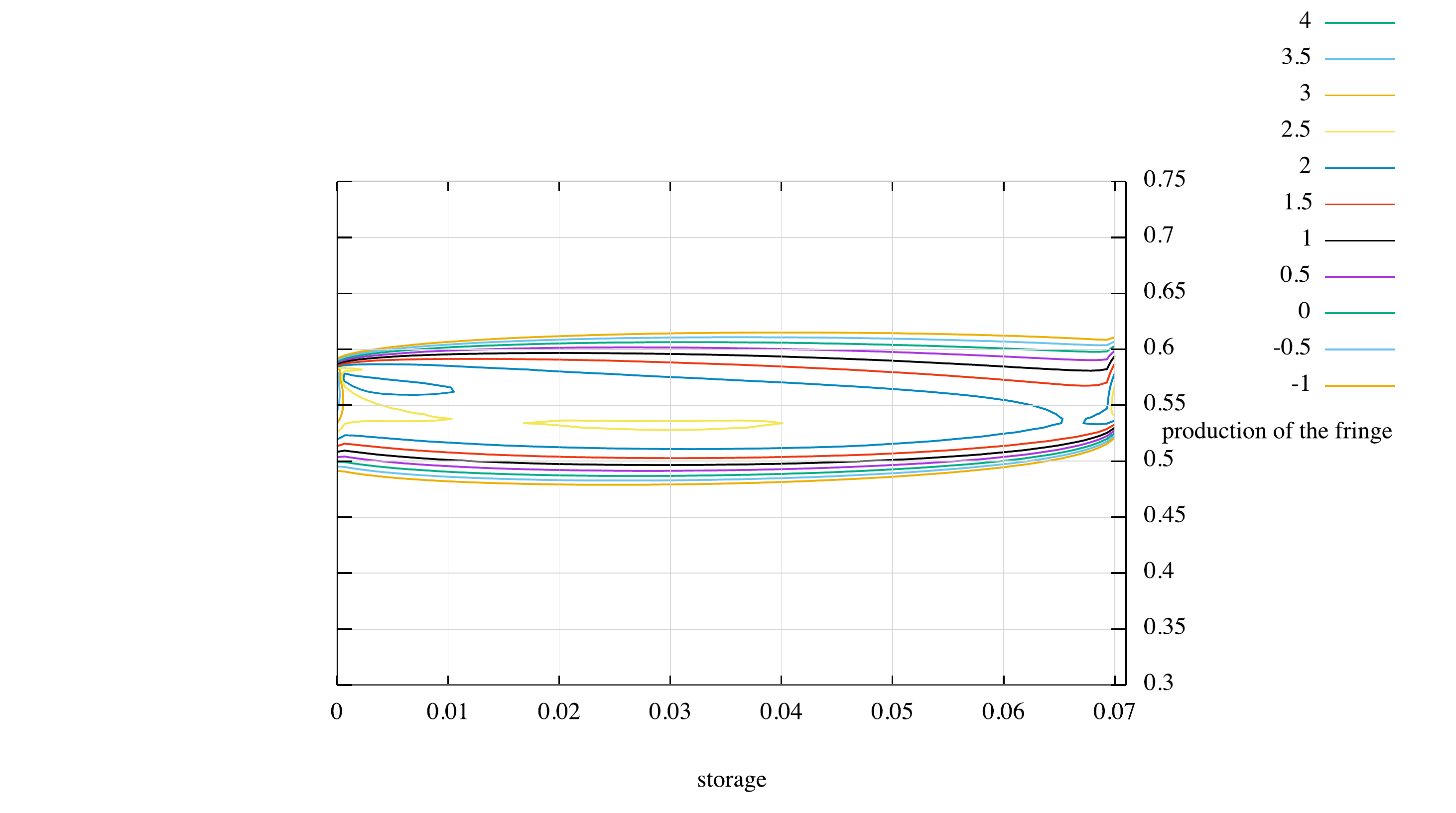}}
  \captionof{figure}{Case 2: the contours of the invariant measure of $(k_t,z_t)$ (in logarithmic scale).}
  \label{fig:25}
\end{center}

\begin{center}
  \includegraphics[width=0.45\linewidth]{{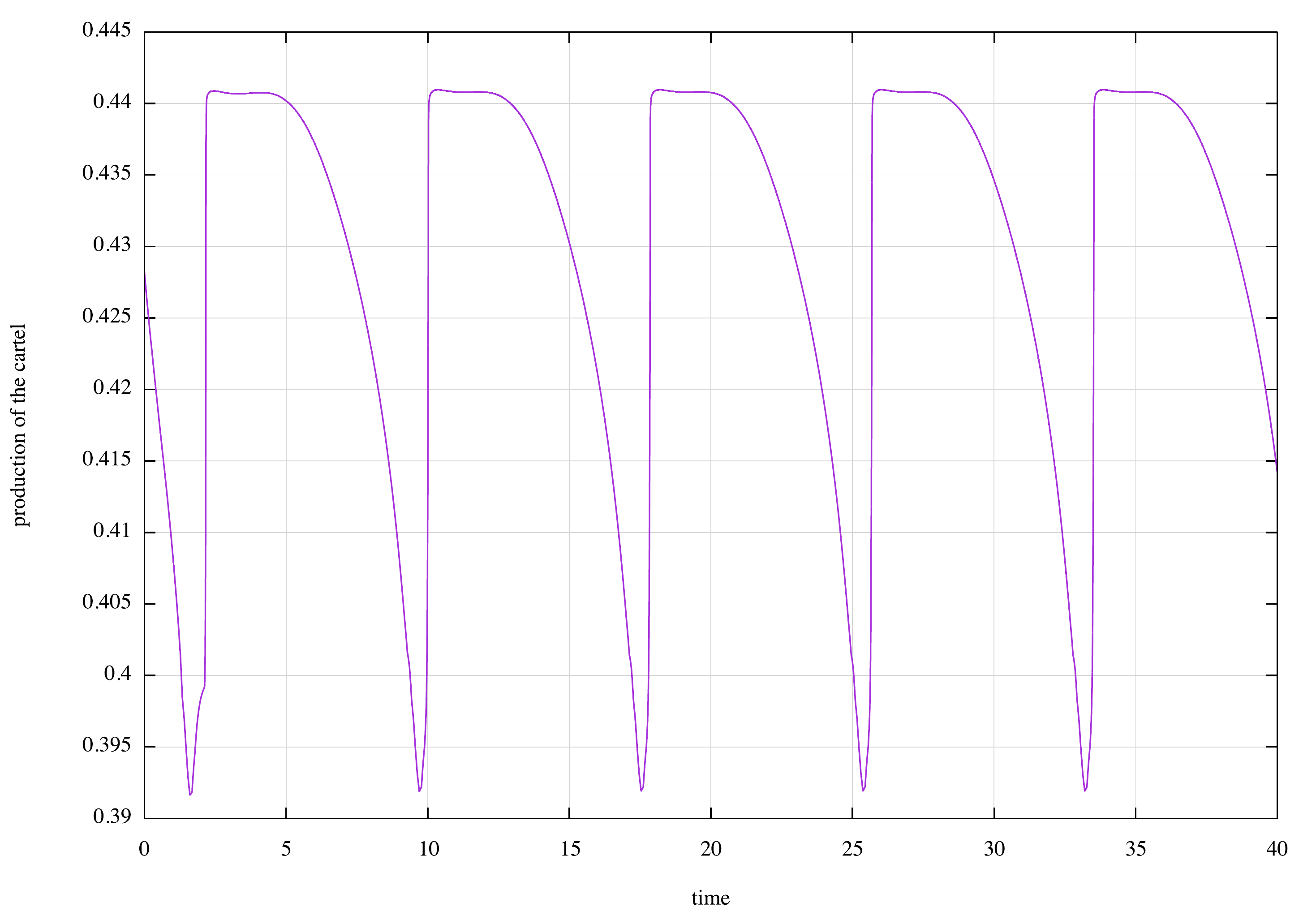}}
\includegraphics[width=0.45\linewidth]{{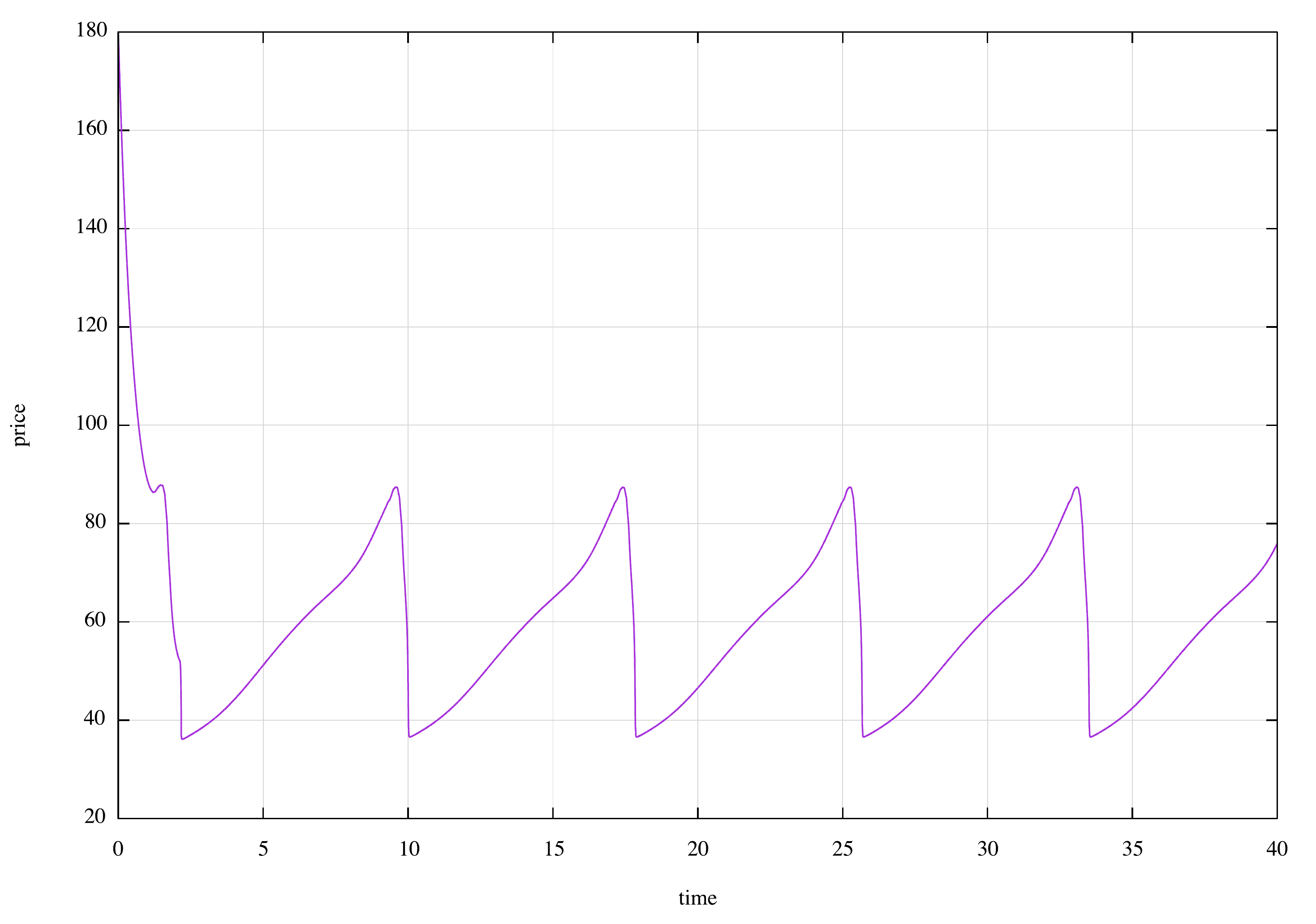}}
\\
\includegraphics[width=0.45\linewidth]{{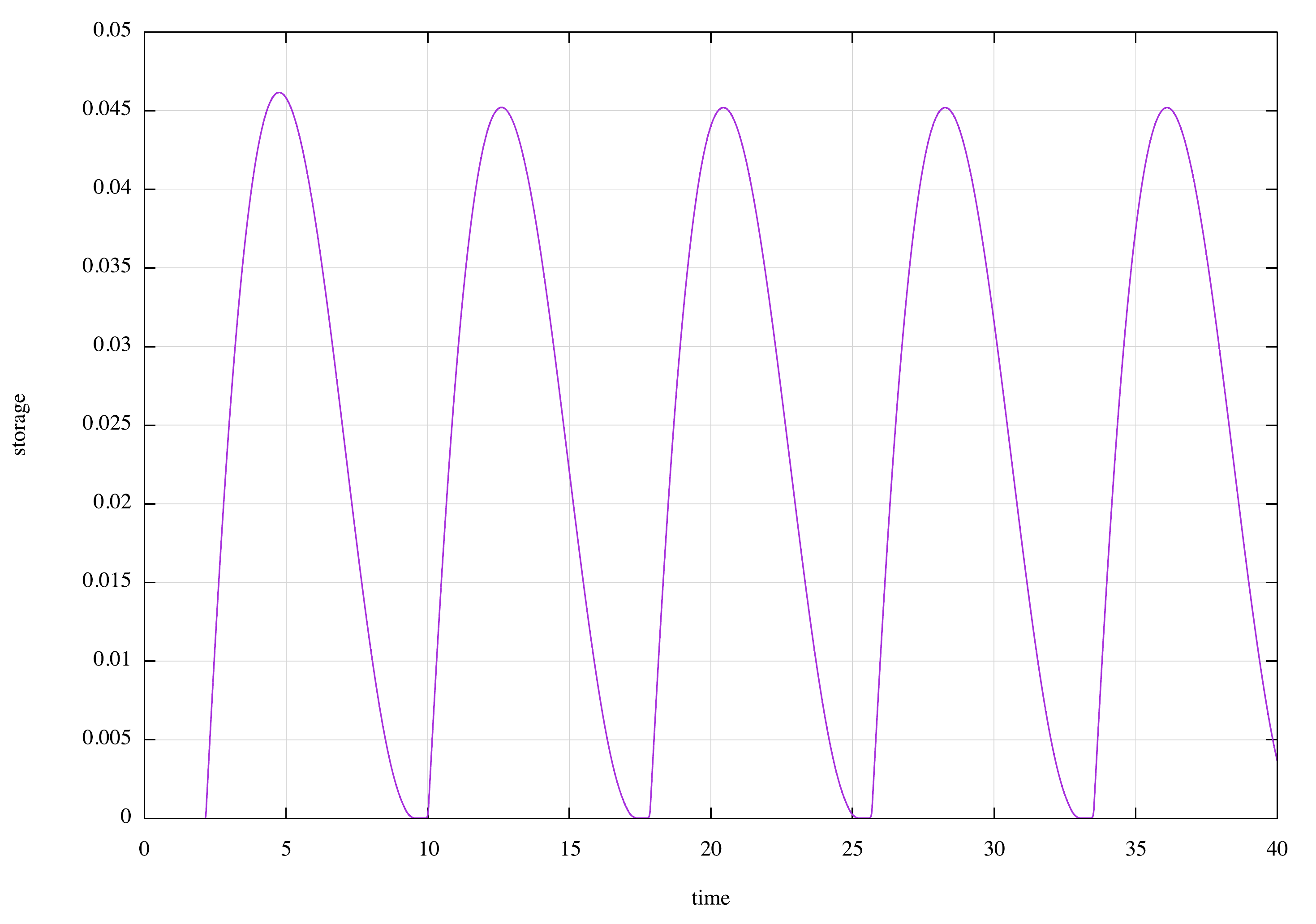}}
\includegraphics[width=0.45\linewidth]{{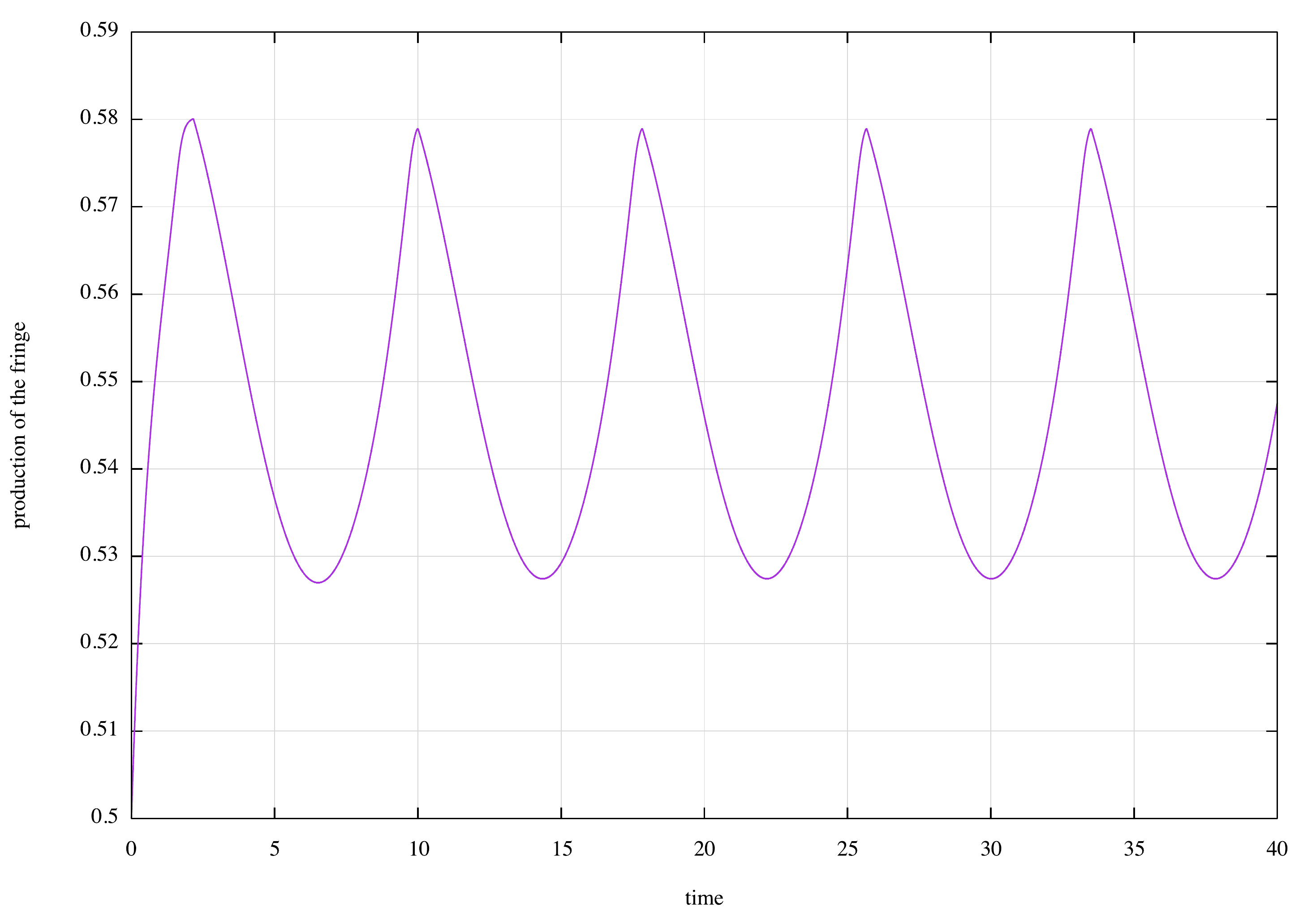}}
  \captionof{figure}{ Case 2 : a stable cycle. Top-Left: the production of the cartel.Top-Right: the price. Bottom-Left: the  level of storage. Bottom-Right: the production of the fringe.}  \label{fig:2cycles}
\end{center}

\end{document}